\documentclass[11pt]{article}
\usepackage{amssymb}
\usepackage{latexsym}
\usepackage{amsfonts}
\usepackage{amsmath}
\newtheorem{thrm}{Theorem}[section]
\newtheorem{lemma}[thrm]{Lemma}

\newtheorem{cor}[thrm]{Corollary}
\newtheorem{remark}[thrm]{Remark}
\newtheorem{exam}{Example}

\numberwithin{equation}{section}
\usepackage[dvips]{color}

\def\P{\mathbb{P} }
\def\Q{\mathbb{Q} }
\def\R{\mathbb{R} }
\def\N{\mathbb{N} }

\def\D{\mathbb{D} }

\def\cA{\mathcal{A} }

\makeatletter
\begin{document}
\allowdisplaybreaks

\title{\Large \bf Williams decomposition for  superprocesses}
\author{ \bf  Yan-Xia Ren\footnote{The research of this author is supported by NSFC (Grant Nos.  11271030 and 11671017)\hspace{1mm} } \hspace{1mm}\hspace{1mm}
Renming Song\thanks{Research supported in part by a grant from the Simons
Foundation (\#429343, Renming Song).} \hspace{1mm}\hspace{1mm} and \hspace{1mm}\hspace{1mm}
Rui Zhang\footnote{The research of this author is supported by NSFC (Grant No. 11601354)}
\hspace{1mm} }
\date{}
\maketitle

\begin{abstract}
We decompose the  genealogy of a general superprocess with spatially dependent
branching mechanism with respect to the last individual alive (Williams decomposition).
This is a generalization of
the main result of Delmas and H\'{e}nard \cite{DH} where only superprocesses with
spatially dependent quadratic
branching mechanism were considered.
As an application of the Williams decomposition, we prove that, for some
superprocesses, the normalized total mass will converge to a point mass at
its extinction time. This generalizes a result of Tribe \cite{Tribe} in the sense that
our branching mechanism is more general.
\end{abstract}

\medskip


\medskip


\section{Introduction}

Let $X$ be a  superprocess with a spatially dependent
branching mechanism.
We assume that the extinction time $H$ of $X$ is finite.
In this paper we study the genealogical structure of $X$. More precisely,
we give a spinal decomposition of $X$ involving the ancestral lineage of the last individual alive, conditioned on $H=h$ with $h>0$ being a constant. This decomposition is called a Williams decomposition,
in analogy with the terminology of
Delmas and H\'{e}nard \cite{DH}.
For a superprocess with spatially independent branching mechanism, the
spatial motion is independent of the genealogical structure.
As a consequence, the
law of the ancestral lineage of the last individual alive
does not depend on the
original motion. Therefore, in this setting, the description of $X$ conditioned on $H=h$ may be deduced
from Abraham and Delmas \cite{AD} where no spatial motion is taken into account.
On the contrary, for a superprocess with nonhomogeneous branching mechanism,
the law of the ancestral lineage
of the last individual alive should depend on the spatial motion and the extinction time $h$.
Delmas and H\'{e}nard \cite{DH} gave a Williams decomposition for
superprocesses with a spatially dependent quadratic branching mechanism given by
$$\Psi(x,z)=\beta(x)z+\alpha(x)z^2,$$
 under some conditions on $\beta(x)$ and $\alpha(x)$ (see ($H2$) and ($H3$) in \cite{DH}).
In \cite{DH}, the Williams decomposition was established for superprocesses
with spatially dependent quadratic branching mechanism by using
two transformations to change the branching mechanism $\Psi(x,z)$
to a spatially independent one,  say $\psi_0$, and  then using the genealogy of
superprocesses with branching mechanism $\psi_0$  given by the Brownian snake.
As mentioned in \cite{DH}, the drawback of the approach in \cite{DH} is that
one has  to restrict  to quadratic branching mechanisms with bounded and smooth
parameters.

The goal of this paper is to establish a Williams decomposition for more general
superprocesses.
Our superprocesses are more general in two aspects: first the spatial motion
can be a general Markov process and secondly the branching mechanism is general and  spatially dependent (see \eqref{e:branm} below).
We will give  conditions
that guarantee our general superprocesses admit a Williams decomposition.
The conditions should be  satisfied by a lot of superprocesses.
We obtain a Williams decomposition by direct construction. For any fixed constant $h>0$, we first describe the motion of a spine up to time $h$ and then construct three  kinds of immigrations (continuous immigration, jump immigration and immigration at time 0) alone the spine. We prove that, conditioned on $H=h$,  the sum
of the contributions of the three types of immigrations
has the same distribution as $X$ before time $h$, see Theorem \ref{the:3.4} below. Note that for quadratic branching mechanisms, there is no jump immigration.

As an application of the Williams decomposition, we prove that, for some
superprocesses, the normalized total mass will converge to a point mass at
its extinction time, see Theorem \ref{mthem} below.
This generalizes a result of Tribe \cite{Tribe} in the sense that
our branching mechanism is more general.

\section{Preliminary}\label{Sec2}
\subsection{Superprocesses and assumptions}

In this subsection,
we describe the superprocesses we are going to work with and
formulate our assumptions.

Suppose that $E$ is a locally compact separable metric
space. Let $E_{\partial}:=E\cup\{\partial\}$ be the one-point compactification
of $E$. $\partial$ will be interpreted as the cemetery point. Any
function $f$ on $E$ is automatically extended to $E_{\partial}$ by setting $f(\partial)=0$.

Let
$\mathbb{D}_E$
be the set of all the c\`{a}dl\`{a}g functions from $[0,\infty)$ into $E_{\partial}$ having $\partial$ as a trap. The filtration is defined by ${\mathcal{F}_t=\mathcal{F}^0_{t+}}$, where ${\mathcal{F}^0_{t}}$ is the natural canonical filtration,
and $\mathcal{F}=\bigvee_{t\ge0}\mathcal{F}_t$.
Consider the canonical process $\xi_t$ on $(\mathbb{D}_E,\{\mathcal{F}_t\}_{t\ge 0})$. We will assume that $\xi=\{\xi_t,\Pi_x\}$ is a Hunt process on $E$ and $\zeta:=
 \inf\{t>0: \xi_t=\partial\}$ is the lifetime of $\xi$.
We will use $\{P_t:t\geq 0\}$ to denote the semigroup of $\xi$.
We will use $\mathcal{B}_b(E)$ ($\mathcal{B}_b^+(E)$) to denote the
set of (non-negative) bounded Borel functions on $E$.
We will use ${\cal M}_F(E)$ to denote the family of finite measures
on $E$
and ${\cal M}_F(E)^{0}$ to denote the family of non-trivial finite measures on $E$.

Suppose that the branching mechanism is given by
\begin{equation}\label{e:branm}
\Psi(x,z)=-\alpha(x)z+b(x)z^2+\int_{(0,+\infty)}(e^{-z y}-1+z y)n(x,dy),
\quad x\in E, \quad z> 0,
\end{equation}
where $\alpha\in \mathcal{B}_b(E)$, $b\in \mathcal{B}_b^+(E)$ and $n$ is a kernel from $E$ to $(0,\infty)$ satisfying
\begin{equation}\label{n:condition}
  \sup_{x\in E}\int_{(0,+\infty)} (y\wedge y^2) n(x,dy)<\infty.
\end{equation}
Then there exists a constant $K>0$, such that
$$
|\alpha(x)|+b(x)+\int_{(0,+\infty)} (y\wedge y^2) n(x,dy)\le K.
$$

  Let ${\cal M}_F(E)$
be the space of finite measures on $E$,
equipped with the topology of weak convergence. As usual,
$\langle f,\mu\rangle:=\int_E f(x)\mu(dx)$
 and $\|\mu\|:=\langle 1,\mu\rangle$.
 According to \cite[Theorem 5.12]{Li11}, there is a Hunt process $X=\{\Omega, {\cal G}, {\cal G}_t, X_t, \P_\mu\}$ taking values in  $\mathcal{M}_F(E)$,
  such that, for every
$f\in \mathcal{B}^+_b(E)$ and $\mu \in \mathcal{M}_F(E)$,
\begin{equation}\label{int}
  -\log \P_\mu\left(e^{-\langle f,X_t\rangle}\right)=\langle u_f(t,\cdot),\mu\rangle,
\end{equation}
where $u_f(t,x)$ is the unique positive solution to the equation
\begin{equation}
   u_f(t,x)+\Pi_x\int_0^t\Psi(\xi_s, u_f(t-s,\xi_s))ds=\Pi_x f(\xi_t),
\end{equation}
where $\Psi(\partial,z)=0$, $z>0$.
 $X=\{X_t:t\ge 0\}$ is called a superprocess
with  spatial motion $\xi=\{\xi_t, \Pi_x\}$ and  branching mechanism
$\Psi$, or sometimes a $(\Psi, \xi)$-superprocess.
In this paper, the superprocess we deal with is always this Hunt realization. For the existence of $X$, see also
\cite{Dawson} and \cite{E.B.}.

Define $v(t,x):=-\log\P_{\delta_x}(\|X_t\|=0)$, and
$H:=\inf\{t\ge 0: \|X_t\|=0\}$.
 It is obvious that $v(0,x)=\infty$.
In this paper, we will consider the critical and subcritical case.
More precisely, throughout this paper, we assume that $X$ satisfy the following uniform global extinction  property.

\begin{itemize}
 \item [{\bf(H1)}]  For any $t>0$,
\begin{equation}\label{4.2'}
 \sup_{x\in E} v(t,x)<\infty\quad\mbox{ and }\quad
 \lim_{t\to\infty}v(t,x)=0.
\end{equation}
\end{itemize}

\begin{remark}
Note that Assumption {\bf(H1)} is equivalent to
\begin{equation}\label{equiv-assum4}\inf_{x\in E}\P_{\delta_x}(\|X_t\|=0)>0\quad \mbox{ for all } t>0\quad\mbox{ and } \quad
\P_{\delta_x}(H<\infty)=\lim_{t\to\infty}\P_{\delta_x}(\|X_t\|=0)=1.
\end{equation}
\end{remark}

\begin{remark}\label{psi-domi} If
\begin{equation}\label{com1}
  \Psi(x,z)\ge \widetilde{\Psi}(z):=bz^2+\int_{0}^\infty \big(e^{-yz}-1+yz\big)n(dy),
\end{equation}
where $b\ge 0$,  $\int_{0}^\infty (y\wedge y^2)n(dy)<\infty$ and  $\widetilde{\Psi}$ satisfies the Grey condition:
$$
\int^\infty \frac{1}{\widetilde{\Psi}(z)}\,dz<\infty,
$$
then  Assumption {\bf(H1)} holds.
\end{remark}

We also assume that

\begin{itemize}
  \item [{\bf(H2)}]For any $x\in E$ and $t>0$,
  \begin{equation}\label{11.1}
   w(t,x):=-\frac{\partial v}{\partial t}(t,x)
\end{equation}
exists.
Moreover, for any $0<r<t$,
\begin{equation}\label{cond3}
  \sup_{r\le s\le t}\sup_{x\in E}w(s,x)<\infty.
\end{equation}
\end{itemize}

Note that, since $t\to v(t,x)$ is decreasing, we have $w(t,x)\ge 0$. We also use $v_t$ and $ w_t$ to denote the function $x\to v(t,x)$ and $x\to w(t,x)$ respectively.

\begin{exam}\label{example1}
Assume that the spatial motion $\xi$ is conservative, that is $P_t(1)\equiv1$,
and the branching mechanism is spatially independent,
that is
$$
\Psi(x,z)=\Psi(z)=az+bz^2+\int_{0}^\infty \big(e^{-yz}-1+yz\big)n(dy),
$$
where $a\ge 0$, $b\ge 0$ and
$\int_{0}^\infty (y\wedge y^2)n(dy)<\infty$. We also assume that $\Psi$ satisfies the Grey condition:
$$
\int^\infty \frac{1}{\Psi(z)}\,dz<\infty.
$$
Then $\{\|X_t\|,t\ge 0\}$ is a continuous state branching process with branching mechanism $\Psi(z)$.
 So $v(t,x)=v(t)<\infty$ does not depend on $x$,
and $\lim_{t\to\infty}v(t)=0$,
thus Assumption {\bf(H1)} holds immediately.
Moreover, for $t>0$, we have that
$$
w(t):=-\frac{\mbox{d}}{\mbox{d}t}v(t)=\Psi(v(t)).
$$
Thus Assumption {\bf(H2)} is satisfied.
See \cite[Theorem 10.1]{Kyprianou} for more details.
\end{exam}

In Section \ref{sec:ex}, we will give more examples, including some class of superdiffusions,
 that satisfy  Assumptions {\bf(H1)}-{\bf(H2)}.

\subsection
{Excursion law of $\{X_t,t\ge 0\}$}

We use $\mathbb{D}$ to denote the space of $\mathcal{M}_F({E})$-valued
c\`{a}dl\`{a}g functions $t\mapsto \omega_t$ on $(0, \infty)$ having zero as a trap.
We use $(\mathcal{A},\mathcal{A}_t)$ to denote the natural $\sigma$-algebras on
$\mathbb{D}$ generated by the coordinate process.

Let $\{Q_{t}(\mu,\cdot):=\P_{\mu}\left(X_{t}\in\cdot\right):t\ge 0,\ \mu\in{\cal M}_F(E)\}$ be the transition
semigroup of $X$. Then by \eqref{int}, we have
\begin{equation}
\int_{{\cal M}_F(E)}e^{-\langle f,\nu\rangle}Q_{t}(\mu,d\nu)=\exp\{-\langle V_{t}f,\mu\rangle\}\quad\mbox{ for }\mu\in{\cal M}_F(E)\mbox{ and }t\ge 0,\nonumber
\end{equation}
where $V_tf(x):=u_f(t,x), x\in E$. This implies that $Q_{t}(\mu_{1}+\mu_{2},\cdot)=Q_{t}(\mu_{1},\cdot)*Q_{t}(\mu_{2},\cdot)$ for any $\mu_{1},\mu_{2}\in{\cal M}_F(E)$, and hence $Q_{t}(\mu,\cdot)$ is an infinitely divisible probability measure on ${\cal M}_F(E)$. By the semigroup property of $Q_{t}$, $V_{t}$ satisfies that
\begin{equation}
V_{s}V_{t}=V_{t+s}\quad\mbox{ for all }s,t\ge 0.\nonumber
\end{equation}
Moreover, by the infinite divisibility of $Q_{t}$, each operator $V_{t}$ has the representation
\begin{equation}
V_{t}f(x)=\lambda_{t}(x,f)+\int_{{\cal M}_F(E)^{0}} \left(1-e^{-\langle f,\nu\rangle}\right)L_{t}(x,d\nu)\quad\mbox{ for } t>0,\ f\in\mathcal{B}^{+}_{b}(E),\label{3.1}
\end{equation}
where $\lambda_{t}(x,dy)$ is a bounded kernel on $E$ and $(1\wedge \nu(1))L_{t}(x,d\nu)$ is a bounded kernel from $E$ to ${\cal M}_F(E)^{0}$.
Let $Q^{0}_{t}$ be the restriction of $Q_{t}$ to
${\cal M}_F(E)^{0}$.
Let $E_{0}:=\{x\in E:\lambda_{t}(x,E)=0\mbox{ for all }t >0\}$.

For $\lambda>0$, we use $V_{t}\lambda$ to denote $V_{t}f$ when the function $f\equiv \lambda$. It then follows from \eqref{3.1} that for every $x\in E$ and $t>0$,
\begin{equation}\nonumber
V_{t}\lambda(x)=\lambda_{t}(x,E)\lambda+\int_{{\cal M}_F(E)^{0}}\left(1-e^{-\lambda\langle 1,\nu\rangle}\right)L_{t}(x,d\nu).
\end{equation}
The left hand side tends to $-\log\P_{\delta_{x}}\left(X_{t}=0\right)$ as $\lambda\to+\infty$.
Therefore, Assumption {\bf(H1)} implies that
 $\lambda_{t}(x,E)=0$ for all $t>0$ and hence $x\in E_{0}$, which says that $E=E_0$.

For $x\in E$,  we get from \eqref{3.1} that
\begin{equation}\nonumber
V_{t}f(x)=\int_{{\cal M}_F(E)^{0}} \left(1-e^{-\langle f,\nu\rangle}\right)L_{t}(x,d\nu)\quad\mbox{ for } t>0,\ f\in\mathcal{B}^{+}_{b}(E).\label{4}
\end{equation}
It then follows from \cite[Proposition 2.8 and Theorem A.40]{Li11} that for every $x\in E$, the family of measures $\{L_{t}(x,\cdot):t>0\}$ on ${\cal M}_F(E)^{0}$ constitutes an entrance law for the restricted semigroup $\{Q^{0}_{t}:t\ge 0\}$.
It is known (see \cite[Section 8.4]{Li11}) that
one can associate with $\{\P_{\delta_x}:x\in E\}$ a family
of $\sigma$-finite measures $\{\N_x:x\in E\}$ defined on $(\mathbb{D},\mathcal{A})$ such  that $\mathbb{N}_x(\{0\})=0$,
\begin{equation}
 \int_{\mathbb{D}}(1- e^{-\langle f, \omega_{t} \rangle})\mathbb{N}_x(d\omega)
= -\log \mathbb{P}_{\delta_x}(e^{-\langle f, X_{t} \rangle}) ,
\quad f\in {\cal B}^+_b(E),\ t> 0,
\label{DK}
\end{equation}
and, for every $0<t_1<\cdots<t_n<\infty$, and nonzero $\mu_1,\cdots,\mu_n\in M_F(E)$,
\begin{eqnarray}\label{TN}
  &&\mathbb{N}_x(\omega_{t_1}\in d\mu_1,\cdots,\omega_{t_n}\in d\mu_n) \nonumber\\
 & =& \mathbb{N}_x(\omega_{t_1}\in d\mu_1)\P_{\mu_1}(X_{t_2-t_1}\in d\mu_2)\cdots \P_{\mu_{n-1}}(X_{t_n-t_{n-1}}\in d\mu_n).
\end{eqnarray}
This measure $\mathbb{N}_{x}$ is called the \textit{Kuznetsov measure} corresponding to the entrance law $\{L_{t}(x,\cdot):t>0\}$ or the \textit{excursion law} for superprocess $X$. For earlier work on
excursion law of superprocesses,
see \cite{E.B2., elk-roe, Li03}.

It follows from \eqref{DK} that for any $t>0$,
\begin{equation}\label{N3}
  \mathbb{N}_x(\|\omega_t\|\neq 0)=-\log\mathbb{P}_{\delta_x}(\|X_t\|=0)<\infty.
\end{equation}

\section{Main results}

In this and the next section we will always assume that
Assumptions {\bf{(H1)}}-{\bf{(H2)}} hold.

Recall that $H:=\inf\{t\ge 0: \|X_t\|=0\}$.
Note that
\begin{equation}\label{def:F}
  F_H(t):=\P_{\mu}(H\le t)=\P_{\mu}(\|X_t\|=0)=e^{-\langle v_t,\mu\rangle}.
\end{equation}
By the continuity of $v(t, x)$ with respect to $t\in(0,\infty)$, we get that for any $t>0$,
\begin{equation}\label{def:F'}\P_{\mu}(H<t)=\lim_{\epsilon\downarrow 0}\P_{\mu}(H\le t-\epsilon)=\lim_{\epsilon\downarrow 0}e^{-\langle v_{t-\epsilon},\mu\rangle}=e^{-\langle v_t,\mu\rangle}=\P_{\mu}(H\le t).\end{equation}

For $h>0$, define
\begin{equation}\label{Martingale1}
  M_t^{h}:=\frac
  {\langle w_{h-t},X_t\rangle e^{-\langle v_{h-t},X_t\rangle}}
  {\langle w_{h},X_0\rangle e^{-\langle v_{h},X_0\rangle}}, \qquad 0\le t<h.
\end{equation}
Then, under $\P_{\mu}$, $\{M_t^h,0\le t<h\}$ is a nonnegative martingale
 with mean one
(see Lemma \ref{lemma10} below).

\begin{thrm}\label{therm1}
For any $h>0$ and $t<h$,
$$\lim_{\epsilon\downarrow0}\P_{\mu} (A|h\le H<h+\epsilon)=
\P_{\mu}(\textbf{1}_AM_t^h),\qquad \forall A\in\mathcal{G}_t.$$
\end{thrm}

We define, for each $h>0$,
$$
\P_{\mu}(\cdot|H=h):=\lim_{\epsilon\downarrow0}\P_{\mu} (\cdot|h\le H<h+\epsilon).
$$
Then, by Theorem  \ref{therm1},
$\{X_t,t<h;\P_{\mu}(\cdot|H=h)\}$ has the same law as $\{X_t,t<h;\P_{\mu}^h\}$, where $\P_{\mu}^h$
is a new measure defined via the martingale $M^h_t$:
$$
\frac{d\P_{\mu}^h}{d\P_{\mu}}\Big|_{\mathcal{G}_t}=M_t^h,\qquad t<h.
$$

\begin{cor}\label{cor1}
For any $A\in \mathcal{G}_t$, we have
$$\P_{\mu}(A\cap\{H>t\})=\int_t^\infty \P_{\mu}^h(A)F_H(dh).$$
\end{cor}
\textbf{Proof:}
It follows from Fubini's theorem that
\begin{align*}
\int_t^\infty \P_{\mu}^h(A)F_H(dh)&=
\int_t^\infty \P_{\mu}(\textbf{1}_AM_t^h)F_H(dh)\\
&=\int_t^\infty \P_{\mu}(\textbf{1}_A\langle w_{h-t},X_t\rangle e^{-\langle v_{h-t},X_t\rangle})\,dh\\
   &= \P_{\mu}\left(\textbf{1}_A\int_t^\infty\langle w_{h-t},X_t\rangle e^{-\langle v_{h-t},X_t\rangle}\,dh\right)\\
   &=\P_{\mu}\left(\textbf{1}_A\int_0^\infty\langle w_{h},X_t\rangle e^{-\langle v_{h},X_t\rangle}\,dh\right)\\
   &=\P_{\mu}\left(A\cap\{X_t\neq0\}\right)=\P_{\mu}\left(A\cap\{H>t\}\right),
\end{align*}
where in the fifth equality  we use the fact that
$$\int_0^\infty\langle w_{h},X_t\rangle e^{-\langle v_{h},X_t\rangle}\,dh=\lim_{h\to\infty}e^{-\langle v_{h},X_t\rangle}-\lim_{h\to0}e^{-\langle v_{h},X_t\rangle}=\textbf{1}_{\{X_t\neq0\}}.$$
\hfill$\Box$

\smallskip

For any $h>0$ and $t\in[0,h)$, we define
$$
Y^{h}_t:=\frac{w(h-t,\xi_t)}{w(h,\xi_0)}e^{-\int_0^t \Psi'_z(\xi_u,v(h-u,\xi_u))\,du},
$$
where $\Psi'_z(x,z)=\frac{\partial \Psi(x,z)}{\partial z}$.
Then we have the following result whose proof will be given in Section \ref{proof}.

\begin{lemma}\label{lemma9}
 Under $\Pi_x$, $\{Y^{h}_t,t<h\}$ is a nonnegative martingale satisfying
$\Pi_x( Y_t^h)= 1$.
\end{lemma}

\begin{remark}
In Example \ref{example1}, $w(t,x)$ and $v(t,x)$ do not depend on $x$, and for any $h>0$ and $0\le t<h$, $Y_t^{h}\equiv1$.
\end{remark}

Now we state our main result: the Williams decomposition. We will
construct a new process
$\{\Lambda_t^h,t<h\}$
which has the same law
as $\{X_t,t<h;\P_{\mu}(\cdot|H=h)\}$.

Let $\mathcal{F}_{h-}:=\bigvee_{t<h}\mathcal{F}_{t}$. Now we define a new probability measure $\Pi_{x}^h$ on $(\mathbb{D}_E,\mathcal{F}_{h-})$ by
$$
\frac{\Pi_{x}^h}{\Pi_{x}}\Big|_{\mathcal{F}_t}:=Y^{h}_t, \qquad t\in [0, h).
$$
Under $\Pi_x^h$, $(\xi_t)_{0\le t<h}$ is a conservative Markov process.
If $\nu$ is a probability measure on $E$, we define
$$
\Pi_{\nu}^h:=\int_E \Pi_{x}^h\,\nu(dx).
$$
Then, under $\Pi_\nu^h$, $(\xi_t)_{0\le t<h}$ is a Markov process with initial
measure $\nu$.

We put
$$
H(\omega):=\inf\{t>0: \|\omega_t\|=0\}, \quad \omega\in\D.
$$

Let $\xi^{h}:=\{(\xi_t)_{0\le t<h},\Pi_\nu^h\}$,
where $\nu(dx)=\frac{w(h,x)}{\langle w(h,\cdot),\mu\rangle}\mu(dx)$.
Given the trajectory of $\xi^{h}$, we  define three processes as follows:
\begin{description}
  \item[Continuous immigration]
Suppose that
$\mathcal{N}^{1,h}(ds,d\omega)$ is a Poisson random measure on $[0,h)\times \mathbb{D}$ with intensity measure
  $2\textbf{1}_{[0,h)}(s) \textbf{1}_{H(\omega)<h-s}b(\xi_s)\mathbb{N}_{\xi_s}(d\omega)ds$. Define, for $t\in[0,h)$,
  \begin{equation}\label{def:X1}
    X^{1,h,\mathbb{N}}_t:=\int_{[0,t)}\int_{\mathbb{D}}\omega_{t-s}\mathcal{N}^{1,h}(ds,d\omega).
  \end{equation}
  \item[Jump immigration] Suppose that $\mathcal{N}^{2,h}(ds,d\omega)$ is a Poisson random measure on $[0,h)\times \mathbb{D}$ with intensity measure
  $\textbf{1}_{[0,h)}(s) \textbf{1}_{H(\omega)<h-s}\int_0^\infty yn(\xi_s,dy)\P_{y\delta_{\xi_s}}(X\in d\omega)\,ds$. Define, for $t\in[0,h)$,
  \begin{equation}\label{def:X2}
    X^{2,h,\P}_t:=\int_{0}^t\int_{\mathbb{D}}\omega_{t-s}\mathcal{N}^{2,h}(ds,d\omega).
  \end{equation}
   \item[Immigration at time 0]\quad Let $\{X^{0,h}_t, 0\le t<h\}$ be a
   process
   distributed according to  the law $\P_{\mu}(X\in \cdot|H<h)$.
\end{description}
We assume that the three processes
$X^{0,h}$, $X^{1,h,\N}$ and $X^{2,h,\P}$ are independent given the trajectory of $\xi^{h}$.
Define
\begin{equation}\label{def:lambda}
  \Lambda_t^h:=X^{0,h}_t+X^{1,h,\N}_t+X^{2,h,\P}_t.
\end{equation}
We write the law of $\Lambda^h$ as $\textbf{P}_{\mu}^{(h)}$.

\begin{thrm}\label{the:3.4}
Under $\textbf{P}_{\mu}^{(h)}$, the process $ \{\Lambda_t^h,t<h\}$ has the same law as $\{X_t,t<h\}$
conditioned on ${H=h}$.
\end{thrm}

If we define $\Lambda^h_t=0,$ for any $t\ge h$, then we get the following result.
\begin{cor}\label{cor2}
$\{X_t;\P_{\mu}\}$ has the same finite dimensional distribution as $$\int_0^\infty \mathbf{P}_{\mu}^{(h)}(\Lambda^h\in\cdot)F_H(dh).$$
\end{cor}
\textbf{Proof:}
Let $f_k\in \mathcal{B}_b^+(E)$, $k=1,2,\cdots,n$ and $0=t_0<t_1<t_2<\cdots<t_n$. We put $t_{n+1}=\infty$
and define $(t_n, t_{n+1}]:=(t_n, \infty)$.
We will show that
$$
  \P_\mu\left(\exp\left\{-\sum_{j=1}^n\langle f_j,X_{t_j}\rangle\right\}\right)=\int_{(0,\infty)} \textbf{P}_{\mu}^{(h)}\left(\exp\left\{-\sum_{j=1}^n\langle f_j,\Lambda_{t_j}^h\rangle\right\}\right)F_H(dh).
$$
Since $\Lambda_t^h=0$, for $t\ge h$,  we get that
\begin{align*}
  &\int_{(0,\infty)} \textbf{P}_{\mu}^{(h)}\left(\exp\left\{-\sum_{j=1}^n\langle f_j,\Lambda_{t_j}^h\rangle\right\}\right)F_H(dh)\\
  =&\sum_{r=0}^{n}
  \int_{(t_r, t_{r+1}]}
   \textbf{P}_{\mu}^{(h)}\left(\exp\left\{-\sum_{j=1}^r\langle f_j,\Lambda_{t_j}^h\rangle\right\}\right)F_H(dh)\\
  =&\sum_{r=0}^{n}
\int_{(t_r, t_{r+1}]}
  \P_{\mu}^{h}\left(\exp\left\{-\sum_{j=1}^r\langle f_j,X_{t_j}\rangle\right\}\right)F_H(dh)\\
  =&\sum_{r=0}^{n} \P_{\mu}\left(\exp\left\{-\sum_{j=1}^r\langle f_j,X_{t_j}\rangle\right\};
  t_r<H\le t_{r+1}\right)\\
  =&\P_\mu\left(\exp\left\{-\sum_{j=1}^n\langle f_j,X_{t_j}\rangle\right\}\right),
\end{align*}
where the second equality follows from  Theorem \ref{the:3.4}, and  the third equality follows from  Corollary \ref{cor1}. The proof is now complete.
\hfill$\Box$

\smallskip

The decomposition \eqref{def:lambda} is called a Williams decomposition or spinal decomposition of the
supperprocess $\{X_t,t<h\}$ conditioned on ${H=h}$,
and $\xi^{h}=\{(\xi_t)_{0\le t<h},\Pi_\nu^h\}$ is called the spine of the decomposition.
It gives us a tool to study the behavior of the superprocesses $X$ near extinction, see  Theorem \ref{mthem} below. To state Theorem \ref{mthem}, we need the following assumption:

\begin{itemize}
  \item [{\bf{(H3)}}]For any bounded open set  $B\subset E$ and any $t>0$, the function
$$
x\to -\log \P_{\delta_x}\Big(\int_0^t X_s(B^c)\,ds =0\Big)
$$
is finite for $x\in B$ and locally bounded.
\end{itemize}

\begin{thrm}\label{mthem}
Assume that {\bf{(H1)}}-{\bf{(H3)}} hold
and that for any $\mu\in\mathcal{M}_F(E)$,
\begin{equation}\label{cond5}
  \lim_{t\uparrow h}\xi_t=\xi_{h-},\quad \Pi_{\nu}^{h}\mbox{-a.s.},
\end{equation}
where $\nu(dx)=\frac{w(h,x)}{\langle w(h,\cdot),\mu\rangle}\mu(dx)$.
Then there exists an  E-valued random variable $Z$ such that
$$
\lim_{t\uparrow H}\frac{X_t}{\|X_t\|}=\delta_Z,\quad \P_{\mu}\mbox{-a.s.},
$$
where the limit above is in the sense of weak convergence. Moreover, conditioned on $\{H=h\}$, $Z$ has the same law as $\{\xi_{h-},\Pi_{\nu}^{h}\}$, that is, for any $f\in\mathcal{C}_b^+(E)$,
\begin{equation}\label{Z:law}
  \P_{\mu}f(Z)=\int_0^\infty \Pi_\nu^h(f(\xi_{h-})) F_H(dh).
\end{equation}
\end{thrm}

Note that, if the martingale
$\{Y_t^{h},0\le t<h\}$
is uniformly integrable, then condition \eqref{cond5} holds. Now we give an example that satisfies
Assumption {\bf{(H3)}}.
\begin{exam}\label{exam6}
Assume that $\xi$ is a diffusion on $\R^d$ with  infinitesimal generator
$$
L=\sum a_{ij}(x)\frac{\partial}{\partial x_i}\frac{\partial}{\partial x_j}+\sum b_j(x)\frac{\partial}{\partial x_j},
$$
which satisfies the following two conditions:
\begin{description}
  \item[(A)] (Uniform ellipticity) There exists a constant $\gamma>0$ such that
$$
 \sum a_{i,j}(x)u_iu_j\ge \gamma\sum u_j^2, \qquad x\in \R^d.
$$
  \item[(B)] $a_{ij}$ and $b_j$ are bounded H\"{o}lder continuous functions.
\end{description}
Suppose that  the branching mechanism $\Psi(x,z)$ satisfies that, for some $\alpha\in(1,2]$ and $c>0$, $\Psi(x,z)\ge c z^{\alpha}$ for all $x\in \R^d$.

Let $\{X,\P_\mu\}$ and $\{\widetilde{X},\widetilde{\P}_\mu\}$ be a $(\xi,\Psi)$-superprocess and a $(\xi, z^\alpha)$-superprocess respectively.
Then, for any open set $B\subset \R^d$,
\begin{align}
  -\log \P_{\delta_x}\left(\exp\left\{-\lambda\int_0^{t}X_s(B^c)\,ds\right\}\right) &\le -\log \widetilde{\P}_{\delta_x}\left(\exp\left\{-\lambda\int_0^{t}\widetilde X_s(B^c)\,ds\right\}\right) \nonumber\\
   &\le -\log \widetilde{\P}_{\delta_x}(\mathcal{R}\subset B),\label{7.23}
\end{align}
where $\mathcal{R}$ is the range of $\widetilde{X}$, which is the minimal closed subset of $\R^d$ which supports all the measures  $\widetilde{X}_t$, $t\ge 0$.
Thus, we have that
$$
-\log \P_{\delta_x}\Big(\int_0^t X_s(B^c)\,ds =0\Big)=\lim_{\lambda\to\infty}-\log \P_{\delta_x}\left(\exp\left\{-\lambda\int_0^{t}X_s(B^c)\,ds\right\}\right)\le -\log \widetilde{\P}_{\delta_x}(\mathcal{R}\subset B).
$$
By \cite[Theorem 8.1]{E.B.},
$x\to -\log \widetilde{\P}_{\delta_x}(\mathcal{R}\subset B)$ is continuous in $x\in B$. Therefore the superprocess $X$
satisfies Assumption {\bf{(H3)}}.
\end{exam}

\begin{remark}
Now we consider the superprocess in Example \ref{example1}. We assume that $\xi$ is a diffusion in $\R^d$ satisfying the conditions in Example \ref{exam6}, and the branching mechanism $\Psi(z)$ satisfies that, for some $\alpha\in(1,2]$ and $c>0$, $\Psi(z)\ge cz^{\alpha}$. Thus, Assumption {\bf{(H3)}} holds.
Since $Y^{h}_t=1$
and $\Pi_{x}^{h}=\Pi_x$,  condition \eqref{cond5} holds automatically.
Therefore, Theorem \ref{mthem} holds  and
$Z$ has the same law as
$\xi_H$, where $\xi_0\sim \nu(dx)=\mu(dx)/\|\mu\|$. Moreover, $\xi$ and $H$ are independent.

Compared with \cite{Tribe}, the example above assumes that the spatial motion $\xi$ is a diffusion, while in \cite{Tribe}, the spatial motion is a Feller process. However, in \cite{Tribe}, the branching mechanism is binary $(\Psi(z)=z^2)$, while in the example above, the branching mechanisms is more general.
\end{remark}

\section{Proofs of Main Results}\label{proof}

 We will use $\P_{r,\delta_x}$ to denote the law of
$X$ starting from the unit mass $\delta_x$ at time $r>0$. Similarly, we will use
$\Pi_{r,x}$ to denote the law of $\xi$ starting from $x$ at time $r>0$. First, we give an useful lemma.

\begin{lemma}\label{lemma13}
Suppose that
$f\in \mathcal{B}_b^+(E)$
and $g_i\in \mathcal{B}_b^+(E)$, $i=1,2,\cdots,n$. For any $0<t_1\le t_2\le \cdots\le t_n$ and $0\le r\le t_n$, we have
\begin{align}\label{11.2}
  &\P_{r,\mu}\left(\langle f,X_{t_n}\rangle \exp\Big\{-\sum_{j:t_j\ge r}\langle g_j,X_{t_j}\rangle\Big\}\right)\nonumber\\
 =&\int_E\Pi_{r,x}\left(\exp\left\{-\int_r^{t_n} \Psi'_z(\xi_u,U_g(u,\xi_u))\,du\right\} f(\xi_{t_n})\right)\,\mu(dx)e^{-\langle U_g(r,\cdot), \mu\rangle},
\end{align}
where
$$
U_g(r,x):=-\log \P_{r,\delta_x}\Big( \exp\Big\{-\sum_{j:t_j\ge r}\langle g_j,X_{t_j}\rangle\Big\}\Big).
$$
In particular, for any $f\in \mathcal{B}_b^+(E)$
and $g\in \mathcal{B}_b^+(E)$, we have
\begin{equation}\label{11.8}
  \P_{\delta_x}\left(\langle f,X_t\rangle e^{-\langle g,X_t\rangle}\right)=\Pi_x\left(\exp\left\{-\int_0^t \Psi'_z(\xi_u,u_g(t-u,\xi_u))\,du\right\} f(\xi_t)\right)e^{-u_g(t,x)}.
\end{equation}
\end{lemma}
\textbf{Proof:}
By \cite[Propersition 5.14]{Li11}, we have that, for $0\le r\le t_n$,
$$
-\log \P_{r,\mu}\left( \exp\Big\{-\sum_{j:t_j\ge r}\langle g_j,X_{t_j}\rangle-\theta\langle f,X_{t_n}\rangle\Big\}\right)=\langle F_{\theta}(r,\cdot),\mu\rangle,
$$
where $F_\theta(r, x)$ is the
unique
bounded positive solution on $[0, t_n]\times E$ of
\begin{equation}\label{11.7}
  F_\theta(r,x)+\Pi_{r,x}\int_r^{t_n} \Psi(\xi_u, F_{\theta}(u,\xi_u))\,du=\sum_{j:t_j\ge r}\Pi_{r,x}g_j(\xi_{t_j})+\theta\Pi_{r,x} f(\xi_{t_n}).
\end{equation}
Let $F'_{\theta}(r,x):=\frac{\partial}{\partial \theta}F_{\theta}(r,x)$.
Then,
$$
\P_{r, \mu}\left(\langle f,X_{t_n}\rangle \exp\Big\{-\sum_{j:t_j\ge r}\langle g_j,X_{t_j}\rangle\Big\}\right)=-\frac{\partial}{\partial \theta}e^{- \langle F_{\theta}(r,\cdot),\mu\rangle}\big|_{\theta=0+}=\langle F'_0(r,\cdot),\mu\rangle e^{-\langle U_g(r,\cdot),\,\mu\rangle}.
$$
Differentiating both sides of \eqref{11.7}
 with respect to $\theta$
and then letting $\theta\to0$, we get that
$$
F'_0(r,x)+\Pi_{r,x}\int_r^{t_n} \Psi'_z(\xi_u, U_g(u,\xi_u))F'_0(u,\xi_u)\,du=\Pi_{r, x} f(\xi_{t_n}),
$$
which implies that
$$
F'_0(r,x)=\Pi_{r, x}\Big[e^{-\int_r^{t_n} \Psi'_z(\xi_u, U_g(u,\xi_u))\,du} f(\xi_{t_n})\Big].
$$
Therefore \eqref{11.2} holds.
\hfill$\Box$

\bigskip
Recall that $v(t,x):=-\log \P_{\delta_x}(\|X_t\|=0)$ and
$w(t,x):=-\frac{\partial v}{\partial t}(t,x)\ge 0$. Recall the definition of $M_t^h$ in \eqref{Martingale1}.

\begin{lemma}\label{lemma10}
 Under $\P_\mu$, $\{M^{h}_t,t<h\}$ is a nonnegative martingale with
$\P_\mu( M_t^h)= 1$.
\end{lemma}
\textbf{Proof:} For any $h>0$ and $0\le t<h$, by Assumption {\bf{(H2)}} and the dominated convergence theorem, we get that
\begin{align}
  \P_{\mu} \Big[\langle w_{h-t},X_t\rangle e^{-\langle v_{h-t},X_t\rangle}\Big]& =\frac{\partial }{\partial h}\P_{\mu}  e^{-\langle v_{h-t},X_t\rangle}\nonumber\\
  &= \frac{\partial }{\partial h} e^{-\langle v_{h},\mu\rangle}
  = \langle w_h,\mu\rangle e^{-\langle v_{h},\mu\rangle},\label{3.34}
\end{align}
where in the second equality,
we used the Markov property of $X$.
Thus, it follows that
$\P_\mu( M_t^h)= 1.$

By the Markov property of $X$, we obtain that, for $s<t<h$,
$$\P_{\mu} \Big[\langle w_{h-t},X_t\rangle e^{-\langle v_{h-t},X_t\rangle}\Big|\mathcal{G}_s] =\P_{X_s} \Big[\langle w_{h-t},X_{t-s}\rangle e^{-\langle v_{h-t},X_{t-s}\rangle}\Big]= \langle w_{h-s},X_s\rangle e^{-\langle v_{h-s},X_s\rangle},$$
which implies that, under $\P_\mu$, $\{M^{h}_t,t<h\}$ is a nonnegative martingale. The proof is complete.
\hfill$\Box$

\bigskip

\noindent\textbf{Proof of Theorem \ref{therm1}:}
For any $A\in\mathcal{G}_t$, by the Markov property of $X$,
\begin{eqnarray*}
  \P_{\mu} (A|h\le H<h+\epsilon) &=& \frac{\P_{\mu} (A\cap \{h\le H<h+\epsilon\})}{\P_{\mu} (h\le H<h+\epsilon)}\\
   &=& \frac{\P_{\mu} (\textbf{1}_A\P_{X_t}(h-t\le H<h-t+\epsilon))}{e^{-\langle v_{h+\epsilon},\mu\rangle}-e^{-\langle v_{h},\mu\rangle}} \\
   &=&  \frac{\P_{\mu} \Big(\textbf{1}_A\big(e^{-\langle v_{h-t+\epsilon},X_t\rangle}-e^{-\langle v_{h-t},X_t\rangle}\big)\Big)}
   {e^{-\langle v_{h+\epsilon},\mu\rangle}-e^{-\langle v_{h},\mu\rangle}}.
\end{eqnarray*}
By Assumption {\bf{(H2)}}, we get that
\begin{equation}\label{3.31}
  \lim_{\epsilon\downarrow0}\frac{1}{\epsilon}\Big(e^{-\langle v_{h+\epsilon},\mu\rangle}-e^{-\langle v_{h},\mu\rangle}\Big)=\langle w_{h},\mu\rangle e^{-\langle v_{h},\mu\rangle}
\end{equation}
and
\begin{equation}\label{3.32}
  \lim_{\epsilon\downarrow0}\frac{1}{\epsilon}\Big(e^{-\langle v_{h-t+\epsilon},X_t\rangle}-e^{-\langle v_{h-t},X_t\rangle}\Big)=\langle w_{h-t},X_t\rangle e^{-\langle v_{h-t},X_t\rangle}.
\end{equation}
Note that, for $0<\epsilon<1$,
\begin{align*}
  &\frac{1}{\epsilon}\Big(e^{-\langle v_{h-t+\epsilon},X_t\rangle}-e^{-\langle v_{h-t},X_t\rangle}\Big)\le \frac{1}{\epsilon}\Big(1-\exp\{-\langle v_{h-t}- v_{h-t+\epsilon},X_t\rangle\}\Big)\\
  &\le \frac{\langle v_{h-t}- v_{h-t+\epsilon},X_t\rangle}{\epsilon}
  \le \sup_{h-t\le s\le h-t+1}\sup_{x\in E}w(s,x)\langle1, X_t\rangle.
\end{align*}
Thus, it follows from the dominated convergence theorem that
\begin{equation}\label{3.33}
\lim_{\epsilon\downarrow0}\frac{1}{\epsilon}
\P_{\mu} \Big(\textbf{1}_A\big(e^{-\langle v_{h-t+\epsilon},X_t\rangle}-e^{-\langle v_{h-t},X_t\rangle}\big)\Big)
=\P_{\mu}\big(\textbf{1}_A\langle w_{h-t},X_t\rangle e^{-\langle v_{h-t},X_t\rangle}\big).
\end{equation}
Thus, by \eqref{3.31} and \eqref{3.33}, we have that
$$\lim_{\epsilon\downarrow0}\P_{\mu} (A|h\le H<h+\epsilon)
=\P_{\mu}(\textbf{1}_AM_t^h).$$
The proof is now complete.
\hfill$\Box$
\bigskip

\noindent\textbf{Proof of Lemma \ref{lemma9}:}\quad
By the Markov property of $X$, we get that,
\begin{equation}\label{veq}
  e^{-v(t+s,x)}=\P_{\delta_x}(X_{t+s}=0)=\P_{\delta_x}(\P_{X_t}(X_{s}=0))=\P_{\delta_x}(e^{-\langle v_s,X_t\rangle}),
\end{equation}
which implies that $u_{v_s}(t,x)=v(t+s,x)$.
By \eqref{3.34} with $h=t+s$ and $\mu=\delta_x$, we get that
\begin{align*}
 w(t+s,x)e^{-v(t+s,x)}&=\P_{\delta_x}(\langle w_s,X_t\rangle e^{-\langle v_s,X_t\rangle})\\
  &=\Pi_x\left(\exp\left\{-\int_0^t \Psi'_z(\xi_u,v(t+s-u,\xi_u))\,du\right\} w(s,\xi_t)\right)e^{-v(t+s,x)},
\end{align*}
where in the last equality we used Lemma \ref{lemma13} and the fact that $u_{v_s}(t,x)=v(t+s,x)$.
Thus, it follows immediately that
\begin{equation}\label{eq:w}
  w(t+s,x)= \Pi_x\left(\exp\left\{-\int_0^t \Psi'_z(\xi_u,v(t+s-u,\xi_u))\,du\right\} w(s,\xi_t)\right).
\end{equation}

For $0<s<t$, by the Markov property of $\xi$, we have that
\begin{eqnarray*}
 &&\displaystyle \Pi_x\Big(w(h-t,\xi_t)e^{-\int_0^t \Psi'_z(\xi_u,v(h-u,\xi_u))\,du}|\mathcal{F}_s\Big)\\
 &=&\displaystyle e^{-\int_0^s \Psi'_z(\xi_u,v(h-u,\xi_u))\,du}\Pi_x\Big(w(h-t,\xi_t)e^{-\int_s^t \Psi'_z(\xi_u,v(h-u,\xi_u))\,du}|\mathcal{F}_s\Big) \\
   &=&\displaystyle  e^{-\int_0^s \Psi'_z(\xi_u,v(h-u,\xi_u))\,du}\Pi_{\xi_s}\Big(w(h-t,\xi_{t-s})e^{-\int_0^{t-s} \Psi'_z(\xi_u,v(h-s-u,\xi_u))\,du}\Big)\\
   &=&\displaystyle e^{-\int_0^s \Psi'_z(\xi_u,v(h-u,\xi_u))\,du}w(h-s,\xi_s),
\end{eqnarray*}
where the last equality above follows from \eqref{eq:w}.
The proof is now complete.\hfill$\Box$

\subsection{Williams decomposition}

\textbf{Proof of Theorem \ref{the:3.4}:}\quad
Let $f_k\in \mathcal{B}_b^+(E)$, $k=1,2,\cdots,n$ and $0=t_0<t_1<t_2<\cdots<t_n=t<h$.
We will show that
$$
  \P_\mu^h\left(\exp\left\{-\sum_{j=1}^n\langle f_j,X_{t_j}\rangle\right\}\right)= \textbf{P}_{\mu}^{(h)}\left(\exp\left\{-\sum_{j=1}^n\langle f_j,\Lambda_{t_j}^h\rangle\right\}\right).
$$

By the definition of $\Lambda_t^h$, we have
\begin{align}
  &\textbf{P}_{\mu}^{(h)}\Big(\exp\Big\{-\sum_{j=1}^n\langle f_j,\Lambda_{t_j}^h\rangle\Big\}\Big)\nonumber \\
   &= \int_E\frac{w(h,x)}{\langle w(h,\cdot),\mu\rangle}\mu(dx)\Pi_x^h \Big[\textbf{P}_{\mu}^{(h)}\Big(\exp\Big\{-\sum_{j=1}^n\langle f_j,\Lambda_{t_j}^h\rangle\Big\}|\xi^h\Big)\Big].\label{3.8}
\end{align}
By the construction of $\Lambda_t^h$, we have
\begin{align}
  &\textbf{P}_{\mu}^{(h)}\Big(\exp\Big\{-\sum_{j=1}^n\langle f_j,\Lambda_{t_j}^h\rangle\Big\}|\xi^h\Big) \nonumber \\
 = &\P_\mu\Big(\exp\Big\{-\sum_{j=1}^n\langle f_j,X_{t_j}\rangle\Big\}|H<h\Big)
  \times \textbf{P}_{\mu}^{(h)}\Big(\exp\Big\{-\sum_{j=1}^n\langle f_j,X^{1,h,\N}_{t_j}\rangle\Big\}|\xi^h\Big)\nonumber \\
  &\quad\times \textbf{P}_{\mu}^{(h)}\Big(\exp\Big\{-\sum_{j=1}^n\langle f_j,X^{2,h,\P}_{t_j}\rangle\Big\}|\xi^h\Big)\nonumber \\
  =:&(I)\times(II)\times(III).\label{3.9}
\end{align}
Define, for $s<h$,
\begin{eqnarray}\label{3.10}
     J_s(h,x)&:=&-\log\P_{\delta_{x}}
   \left[e^{-\sum_{j=1}^n\langle f_j,X_{t_j-s}\rangle\textbf{1}_{s\le t_j}};\|X_{h-s}\|=0\right].
\end{eqnarray}
We first deal with part (I). By \eqref{3.10}, we have
\begin{equation}\label{3.4}
  J_0(h,x)=-\log \P_{\delta_x}\Big(\exp\Big\{-\sum_{j=1}^n\langle f_j,X_{t_j}\rangle\Big\};\|X_h\|=0\Big).
\end{equation}
By \eqref{def:F'}, $\P_{\mu}(H<h)=\P_{\mu}(H\le h)=e^{-\langle v_h,\,\mu\rangle}$.
Thus we have
\begin{equation}\label{3.3}
  (I)=e^{\langle v(h,\cdot),\,\mu\rangle}e^{-\langle J_0(h,\cdot),\,\mu\rangle}.
\end{equation}
Next we deal with part (II).  By the definition of $X^{1,h,\N}$
and Fubini's theorem, we have
\begin{align}
  \sum_{j=1}^n\langle f_j,X^{1,h,\N}_{t_j}\rangle
  &=\sum_{j=1}^n \int_0^{t}\int_{\D}\langle f_j,\omega_{t_j-s}\rangle\textbf{1}_{s< t_j}\mathcal{N}^{1,h}(ds,d\omega)\nonumber\\
  &=\int_0^{t}\int_{\D}\sum_{j=1}^n\langle f_j,\omega_{t_j-s}\rangle\textbf{1}_{s< t_j}\,\mathcal{N}^{1,h}(ds,d\omega).\label{3.5}
\end{align}
Therefore,
\begin{eqnarray*}
  (II) &=& \textbf{P}_{\mu}^{(h)}\left(\exp\left\{-\int_0^{t}\int_{\D}\sum_{j=1}^n\langle f_j,\omega_{t_j-s}\rangle\textbf{1}_{s< t_j}\,\mathcal{N}^{1,h}(ds,d\omega)\right\}|\xi^h\right) \\
   &=& \exp\left\{-\int_0^{t} 2b(\xi_s)\,ds\int_{\D}\left(1-e^{-\sum_{j=1}^n\langle f_j,\,\omega_{t_j-s}\rangle\textbf{1}_{s< t_j}}\right)\textbf{1}_{H(\omega)<h-s}\N_{\xi_s}(d\omega)\right\}.
\end{eqnarray*}
By the dominated convergence theorem, we obtain that, for $s\neq t_j$, $j=1,2,\cdots,n$,
\begin{align*}
 &\int_{\D}\left(1-e^{-\sum_{j=1}^n\langle f_j,\omega_{t_j-s}\rangle\textbf{1}_{s< t_j}}\right)\textbf{1}_{H(\omega)<h-s}\N_{\xi_s}(d\omega)\\
   =&\int_{\D}\left(1-e^{-\sum_{j=1}^n\langle f_j,\omega_{t_j-s}\rangle\textbf{1}_{s< t_j}}\right)\textbf{1}_{\|\omega_{h-s}\|=0}\N_{\xi_s}(d\omega) \\
   =&\lim_{\theta\to\infty}\int_{\D}\left(1-e^{-\sum_{j=1}^n\langle f_j,\omega_{t_j-s}\rangle\textbf{1}_{s< t_j}}\right)e^{-\theta\|\omega_{h-s}\|}\N_{\xi_s}(d\omega) \\
   =&\lim_{\theta\to\infty}\int_{\D}\left(1-e^{-\sum_{j=1}^n\langle f_j,\omega_{t_j-s}\rangle\textbf{1}_{s< t_j}-\theta\|\omega_{h-s}\|}\right)\N_{\xi_s}(d\omega)-
   \int_{\D}\left(1-e^{-\theta\|\omega_{h-s}\|}\right)\N_{\xi_s}(d\omega)\\
   =&\lim_{\theta\to\infty}-\log\P_{\delta_{\xi_s}}e^{-\sum_{j=1}^n\langle f_j,X_{t_j-s}\rangle\textbf{1}_{s< t_j}-\theta\|X_{h-s}\|}+\log\P_{\delta_{\xi_s}}e^{-\theta\|X_{h-s}\|}\\
   =&-\log\P_{\delta_{\xi_s}}
   \left[e^{-\sum_{j=1}^n\langle f_j,X_{t_j-s}\rangle\textbf{1}_{s< t_j}};\|X_{h-s}\|=0\right]
   +\log\P_{\delta_{\xi_s}}\left(\|X_{h-s}\|=0\right)\\
   =&J_s(h,\xi_s)-v(h-s,\xi_s).
\end{align*}
Hence,
\begin{equation}\label{3.6}
  (II)=\exp\left\{-\int_0^{t} 2b(\xi_s)\Big(J_s(h,\xi_s)-v(h-s,\xi_s)\Big)\,ds\right\}.
\end{equation}
Now we deal with (III). Using arguments similar to those leading to \eqref{3.5}, we get that
 $$
 \sum_{j=1}^n\langle f_j,X^{2,h,\P}_{t_j}\rangle=\int_0^{t}\int_{\D}\sum_{j=1}^n\langle f_j,\omega_{t_j-s}\rangle\textbf{1}_{s\le t_j}\,\mathcal{N}^{2,h}(ds,d\omega).
 $$
Thus,
\begin{align}\label{3.7}
  (III) &= \textbf{P}_{\mu}^{(h)}\left(\exp\left\{-\int_0^{t}\int_{\D}\sum_{j=1}^n\langle f_j,\omega_{t_j-s}\rangle\textbf{1}_{s\le t_j}\,\mathcal{N}^{2,h}(ds,d\omega)\right\}|\xi^h\right) \nonumber\\
   &= \exp\left\{-\int_0^{t} \,ds\int_0^\infty y n(\xi_s,dy)\P_{y\delta_{\xi_s}}\left[\Big(1-e^{-\sum_{j=1}^n\langle f_j,\,X_{t_j-s}\rangle\textbf{1}_{s\le t_j}}\Big)\textbf{1}_{H<h-s}\right]\right\}\nonumber\\
   &=\exp\left\{-\int_0^{t} \,ds\int_0^\infty y n(\xi_s,dy)\Big(e^{-yv(h-s,\,\xi_s)}-e^{-yJ_s(h,\,\xi_s)}\Big)\right\}.
\end{align}
Recall that
$$
\Psi'_z(x,z)=-\alpha(x)+2b(x)z+\int_0^\infty y(1-e^{-yz})n(x,dy).
$$
Combining \eqref{3.6} and \eqref{3.7}, we get that
\begin{align}\label{3.21}
  &(II)\times(III)\nonumber\\
   =  &\exp\left\{-\int_0^{t}\Big(2b(\xi_s)J_s(h,\xi_s)+\int_0^\infty y \Big(1-e^{-yJ_s(h,\xi_s)}\Big)n(\xi_s,dy)\Big)\,ds\right\}\nonumber\\
   &\quad \times\exp\left\{\int_0^{t}\Big(2b(\xi_s)v(h-s,\xi_s)-\int_0^\infty y \Big(1-e^{-yv(h-s,\xi_s)}\Big)n(\xi_s,dy)\Big)\,ds\right\}\nonumber\\
   =&\exp\left\{-\int_0^{t}\Psi'_z(\xi_s,J_s(h,\xi_s))\,ds\right\}\times \exp\left\{\int_0^{t}\Psi'_z(\xi_s,v(h-s,\xi_s))\,ds\right\}.
\end{align}
By \eqref{3.9}, \eqref{3.3} and \eqref{3.21}, we get that, for $h>t$,
\begin{align*}
&\Pi_x^h\Big[\textbf{P}_{\mu}^{(h)}\left(\exp\left\{-\sum_{j=1}^n \langle f_j,\Lambda_{t_j}^h\rangle\right\}|\xi^h\right)\Big]\\
=&e^{\langle v(h,\cdot),\,\mu\rangle}e^{-\langle J_0(h,\cdot),\,\mu\rangle}\Pi_x^h\Big[\exp\left\{-\int_0^{t}\Psi'_z(\xi_s,J_s(h,\xi_s))\,ds\right\}\times \exp\left\{\int_0^{t}\Psi'_z(\xi_s,v(h-s,\xi_s))\,ds\right\}\Big]\\
=&e^{\langle v(h,\cdot),\mu\rangle}e^{-\langle J_0(h,\cdot),\,\mu\rangle}\Pi_x\Big[\frac{w(h-t,\xi_{t})}{w(h,x)}\exp\left\{-\int_0^t\Psi'_z(\xi_s,J_s(h,\xi_s))\,ds\right\}\Big].
\end{align*}
So, by \eqref{3.8}, we obtain that
\begin{align}
  & \textbf{P}_{\mu}^{(h)}\left(\exp\left\{-\sum_{j=1}^n\langle f_j,\Lambda_{t_j}^h\rangle\right\}\right) \nonumber\\
   =&\frac{e^{\langle v(h,\cdot),\mu\rangle}}{\langle w_h,\mu\rangle}
   e^{-\langle J_0(h,\cdot),\,\mu\rangle}\int_E\Pi_x\Big[w(h-t,\xi_{t})\exp\left\{-\int_0^t\Big(\Psi'_z(\xi_s,J_s(h,\xi_s))\Big)\,ds\right\}\Big] \,\mu(dx).
   \label{3.16}
\end{align}
Now we calculate $J_s(h,x)$ defined in \eqref{3.10}.
For $0\le s<t<h$, by the Markov property of $X$, we have that
\begin{align}
 J_s(h,x)& =
  -\log \P_{\delta_{x}}\left[e^{-\sum_{j=1}^n\langle f_j,X_{t_j-s}\rangle\textbf{1}_{s\le t_j}}\P_{X_{t-s}}(\|X_{h-t}\|=0)\right] \nonumber\\
  &=-\log \P_{\delta_{x}}\left[e^{-\sum_{j=1}^n\langle f_j,X_{t_j-s}\rangle\textbf{1}_{s\le t_j}-\langle v(h-t,\cdot),X_{t-s}\rangle}\right]\nonumber\\
  &=-\log \P_{s,\delta_{x}}\left[e^{-\sum_{j=1}^n\langle f_j,X_{t_j}\rangle\textbf{1}_{s\le t_j}-\langle v(h-t,\cdot),X_{t}\rangle}\right].\label{3.19}
\end{align}
Using Lemma \ref{lemma13} with $r=0$, we have that
\begin{align*}
 &e^{-\langle J_0(h,\cdot),\,\mu\rangle}\int_E\Pi_x\Big[w(h-t,\xi_{t})\exp\left\{-\int_0^t\Big(\Psi'_z(\xi_s,J_s(h,\xi_s))\Big)\,ds\right\}\Big] \,\mu(dx)\\
=&\P_{\mu}\Big[\langle w(h-t,\cdot),X_t\rangle \exp\Big\{-\sum_{j=1}^n\langle f_j,X_{t_j}\rangle-\langle v(h-t,\cdot),X_{t}\rangle\Big\}\Big].
\end{align*}
Thus, by \eqref{3.16}, we get that
\begin{equation*}
   \textbf{P}_{\mu}^{(h)}\Big(\exp\Big\{-\sum_{j=1}^n\langle f_j,\Lambda_{t_j}^h\rangle\Big\}\Big)
   =  \P_{\mu}\Big[\exp\Big\{-\sum_{j=1}^n\langle f_j,X_{t_j}\rangle\Big\} M_t^h\Big].
\end{equation*}
Now, the proof is complete.

\hfill$\Box$

\subsection{The behavior of $X_t$ near extinction}

Recall that, for any $\mu\in \mathcal{M}_F(E)$,
$\xi^{h}=\{(\xi_t)_{0\le t<h},\Pi_\nu^h\}$,
where $\nu(dx)=\frac{w(h,x)}{\langle w(h,\cdot),\mu\rangle}\mu(dx)$.

\begin{lemma} \label{the:4.1}
Suppose that Assumptions
{\bf{(H1)}}-{\bf{(H3)}}) hold and that for any $\mu\in\mathcal{M}_F(E)$,
\begin{equation*}
  \lim_{t\uparrow h}\xi_t=\xi_{h-},\quad \Pi_{\nu}^{h}\mbox{-a.s.},
\end{equation*}
where $\nu(dx)=\frac{w(h,x)}{\langle w(h,\cdot),\mu\rangle}\mu(dx)$. Then, for any $h>0$,
$$
\lim_{t\uparrow h}\frac{\Lambda^h_t}{\|\Lambda^h_t\|}=\delta_{\xi_{h-}},\qquad\textbf{P}_\mu^{(h)}\mbox{-a.s.}
$$
\end{lemma}
\textbf{Proof:} By the decomposition \eqref{def:lambda},  we have
$$ \Lambda_t^h:=X^{0,h}_t+X^{1,h,\N}_t+X^{2,h,\P}_t.$$
Define
$$
H_0:=\inf\{t\ge0:X^{0,h}_t=0\}\quad and \quad H(\Lambda^h):=\inf\{t\ge0:\Lambda^{h}_t=0\}.
$$
Then by the definition of $X^{0,h}$, we have $H_0<h$.
By Theorem \ref{the:3.4}, $H(\Lambda^h)=h$.
It follows that
\begin{equation}\label{w5.1}
  \lim_{t\uparrow h}\frac{X^{0,h}_t}{\|\Lambda^h_t\|}=0,\qquad\textbf{P}_\mu^{(h)}\mbox{-a.s.}
\end{equation}
Note that $E_\partial$ is a compact separable metric space.
According to \cite[Exercise 9.1.16 (iii)]{S}, $C_b(E_\partial;\R)$, the space of bounded
continuous $\R$-valued functions $f$ on $E_\partial$, is separable. Therefore, $C_b^+(E)$,  the space of nonnegative bounded
continuous $\R$-valued functions $f$ on $E$,  is also
a separable space.
It suffices to  prove that, for any $f\in C_b^+(E)$,
\begin{equation}\label{toprove1}
\textbf{P}_\mu^{(h)}\left(\lim_{t\uparrow h}\frac{\langle f_h,X^{1,h,\N}_t\rangle+\langle f_h,X^{2,h,\P}_t\rangle}{\|\Lambda^h_t\|}=0\right)=1,
\end{equation}
where $f_h(x)=f(x)-f(\xi_{h-})$.
Note that
$$\textbf{P}_\mu^{(h)}\left(\lim_{t\uparrow h}\frac{\langle f_h,X^{1,h,\N}_t\rangle+\langle f_h,X^{2,h,\P}_t\rangle}{\|\Lambda^h_t\|}=0\right)
=\textbf{P}_\mu^{(h)}\left[\textbf{P}_\mu^{(h)}\left(\lim_{t\uparrow h}\frac{\langle f_h,X^{1,h,\N}_t\rangle+\langle f_h,X^{2,h,\P}_t\rangle}{\|\Lambda^h_t\|}=0\big|\xi^h\right)\right].
$$
Therefore, it suffices to  prove that, for any $f\in C_b^+(E)$,
\begin{equation}\label{toprove}
\textbf{P}_\mu^{(h)}\left(\lim_{t\uparrow h}\frac{\langle f_h,X^{1,h,\N}_t\rangle+\langle f_h,X^{2,h,\P}_t\rangle}{\|\Lambda^h_t\|}=0\big|\xi^h\right)=1, \quad \textbf{P}_\mu^{(h)}\mbox{-a.s.}
\end{equation}

{\bf Step 1} We first prove that given $\xi^h$,
\begin{equation}\label{step1}
 \lim_{t\uparrow h}\frac{\langle f_h, X_t^{1,h,\N}\rangle}{\|\Lambda_t^{h}\|}=0,\quad \textbf{P}_\mu^{(h)}\mbox{-a.s.}
\end{equation}
Note that given $\xi^h$,
$$
\langle f_h,X^{1,h,\mathbb{N}}_t\rangle:=\int_{0}^t\int_{\mathbb{D}}\langle f_h,\omega_{t-s}\rangle \mathcal{N}^{1,h}(ds,d\omega),
$$
where $\mathcal{N}^{1,h}(ds,d\omega)$ is a Poisson random measure on $[0,h)\times \mathbb{D}$ with intensity measure
$$
 2\textbf{1}_{[0,h)}(s) \textbf{1}_{H(\omega)<h-s}b(\xi_s)\mathbb{N}_{\xi_s}(d\omega)ds.
 $$
Let $I_1$ be the support of the measure $\mathcal{N}^{1,h}$. Note that $I_1$ is a random subset of $[0,h)\times \mathbb{D}$.

In the remainder of this proof, we always assume that $\xi^h$ is given.
Since $f\in C_b^+(E)$, for any $\epsilon>0 $,
there exists $\delta_1>0$, depending on $ \xi_{h-}$,
such that $|f(x)-f(\xi_{h-})|\le \epsilon$ for all $|x-\xi_{h-}|\le \delta_1$.
It follows from the fact that $\xi_{h-}=\lim_{s\uparrow h}\xi_s$
there exists $\delta_2\in(0, h)$, depending on $ \xi_{h-}$,
such that $|\xi_s-\xi_{h-}|< \delta_1/2$ for all $s\in(h-\delta_2,h)$. Let $B:=B(\xi_{h-},\delta_1)=\{x\in E:|x-\xi_{h-}|< \delta_1\}$.
Then, for any $t\in(h-\delta_2/2,h)$, we have
\begin{align}
  |\langle f_h,X^{1,h,\N}_t\rangle|=&|\langle f_h\textbf{1}_{\bar{B}},X^{1,h,\N}_t\rangle+\langle f_h\textbf{1}_{\bar{B}^c},X^{1,h,\N}_t\rangle| \nonumber\\
   \le & \epsilon \langle 1,X^{1,h,\N}_t\rangle+2\|f\|_\infty \langle \textbf{1}_{\bar{B}^c},X^{1,h,\N}_t\rangle \nonumber\\
   \le &\epsilon \langle 1,\Lambda^h_t\rangle+2\|f\|_\infty\int_{0}^{h-\delta_2}\int_{\mathbb{D}}\langle 1,\omega_{t-s}\rangle \mathcal{N}^{1,h}(ds,d\omega)\nonumber\\
   &\quad+2\|f\|_\infty\int_{h-\delta_2}^t\int_{\mathbb{D}}\langle \textbf{1}_{\bar{B}^c},\omega_{t-s}\rangle \mathcal{N}^{1,h}(ds,d\omega)\nonumber\\
   =:&\epsilon \langle 1,\Lambda^h_t\rangle+2\|f\|_\infty J_1(t)+2\|f\|_\infty J_2(t).\label{7.11}
\end{align}
It follows that
\begin{eqnarray}\label{w5.2}
  \frac{|\langle f_h,X^{1,h,\N}_t\rangle|}{\|\Lambda^h_t\|} &\le & \epsilon+2\|f\|_\infty\frac{J_1(t)}{\|\Lambda^h_t\|} +2\|f\|_\infty\frac{J_2(t)}{\|\Lambda^h_t\|}.
\end{eqnarray}

First we deal with  $J_1$. For $s\in(0,h-\delta_2)$ and $t\in(h-\delta_2/2,h)$, we have $t-s>\delta_2/2$.
Thus, for $t\in(h-\delta_2/2,h)$, we have
$$
J_1(t)=\int_{0}^{h-\delta_2}\int_{\omega(\delta_2/2)\neq 0,\,H(\omega)<h-s}\langle 1,\omega_{t-s}\rangle \mathcal{N}^{1,h}(ds,d\omega)=\sum_{(s,\,\omega)\in (I_1\cap S_1)}\langle 1,\omega_{t-s}\rangle,
$$
where
\begin{equation}\label{e:S_1}
S_1:=\{(s,\omega):s\in [0,h-\delta_2), \,w(\delta_2/2)\neq 0\quad  \mbox{and} \quad H(\omega)<h-s\}.
\end{equation}
Notice  that
\begin{align}\label{w5.3}
  &\int_{S_1}2\textbf{1}_{[0,h)}(s) \textbf{1}_{H(\omega)<h-s}b(\xi_s)\mathbb{N}_{\xi_s}(d\omega)ds\nonumber\\
 \le  &2K\int_{0}^{h-\delta_2}\mathbb{N}_{\xi_s}(w(\delta_2/2)\neq 0)ds\nonumber\\
  =&2K\int_{0}^{h-\delta_2}v(\delta_2/2,\xi_s)ds\le 2Kh\|v_{\delta_2/2}\|_\infty<\infty,
\end{align}
which implies that given $\xi^h$,
$$
\mathcal{N}^{1,h}(S_1)<\infty,\quad\textbf{P}_\mu^{(h)}\mbox{-a.s.}
$$
That is, given $\xi^h$, $\sharp\{I_1\cap S_1\}<\infty$, $\textbf{P}_\mu^{(h)}$-a.s.
For any $(s,\omega)\in (I_1\cap S_1)$, we have $s+H(\omega)<h$, which implies  that $H_1:=\max_{(s,\omega)\in (I_1\cap S_1)}(s+H(\omega))<h$.
Thus, for any $t\in (H_1, h)$, $J_1(t)=0$, which implies that given $\xi^h$,
\begin{equation}\label{w5.4}
  \lim_{t\uparrow h}\frac{J_1(t)}{\|\Lambda^h_t\|} =0,\quad\textbf{P}_\mu^{(h)}\mbox{-a.s.}
\end{equation}

To deal with $J_2$,
we define
\begin{equation}\label{S_2}
  \mathbb{D}_1:=\{\omega:\exists u\in(0,\delta_2), \mbox{such that } \langle \textbf{1}_{\bar{B}^c},\omega_u\rangle>0\},\quad \mbox{and } \quad S_2=[h-\delta_2,h)\times \mathbb{D}_1.
\end{equation}
Then,
$$
J_2(t)=\sum_{(s,\,\omega)\in (I_1\cap S_2)}\langle \textbf{1}_{\bar{B}^c},\omega_{t-s}\rangle\textbf{1}_{s<t}\,.
$$
We claim that
$\sharp\{I_1\cap S_2\}<\infty.$
Then using arguments similar to those leading to \eqref{w5.4}, we can get that given $\xi^h$,
\begin{equation}\label{w5.5}
  \lim_{t\uparrow h}\frac{J_2(t)}{\|\Lambda^h_t\|} =0,\quad\textbf{P}_\mu^{(h)}\mbox{-a.s.}
\end{equation}
Now we prove the claim. It suffices to prove that given $\xi^h$
\begin{equation}\label{w5.6}
  \int_{S_2}2\textbf{1}_{[0,h)}(s) \textbf{1}_{H(\omega)<h-s}b(\xi_s)\mathbb{N}_{\xi_s}(d\omega)ds<\infty.
\end{equation}
Note that
$$
\int_{S_2}2\textbf{1}_{[0,h)}(s) \textbf{1}_{H(\omega)<h-s}b(\xi_s)\mathbb{N}_{\xi_s}(d\omega)ds\le 2K\int_{h-\delta_2}^h\mathbb{N}_{\xi_s}(\mathbb{D}_1)ds.
$$
For $\omega\in\mathbb{D}$, we have
\begin{align*}
  \mathbb{D}_1=&\{\omega\in\mathbb{D}:\exists u\in(0,\delta_2), \mbox{such that } \langle \textbf{1}_{\bar{B}^c},\omega_u\rangle>0\}\\
  =&\left\{\omega\in\mathbb{D}:\int_0^{\delta_2} \langle \textbf{1}_{\bar{B}^c},\omega_u\rangle\,du>0\right\}\\
  \subset&\left\{\omega\in\mathbb{D}:\int_0^{\delta_2} \langle \textbf{1}_{B^c},\omega_u\rangle\,du>0\right\}.
\end{align*}
Thus,
\begin{align}
  \mathbb{N}_x(\mathbb{D}_1)
  \le&\mathbb{N}_x\left(\int_0^{\delta_2}\langle \textbf{1}_{B^c},\omega_{u}\rangle\,du>0\right)\nonumber\\
  =&\lim_{\lambda\to\infty}\mathbb{N}_x\left(1-\exp\left\{-\lambda\int_0^{\delta_2}\langle \textbf{1}_{B^c},\omega_{u}\rangle\,du\right\}\right)\nonumber\\
  =&\lim_{\lambda\to\infty}-\log \P_{\delta_x}\left(\exp\left\{-\lambda\int_0^{\delta_2}\langle \textbf{1}_{B^c},X_{u}\rangle\,du\right\}\right)\nonumber\\
  =&-\log \P_{\delta_x}\left(\int_0^{\delta_2}\langle \textbf{1}_{B^c},\omega_{u}\rangle\,du=0\right).\label{7.22}
\end{align}
Combining \eqref{7.22} and Assumption {\bf{(H3)}}, we get
\begin{align*}
  &\int_{S_2}2\textbf{1}_{[0,h)}(s) \textbf{1}_{H(\omega)<h-s}b(\xi_s)\mathbb{N}_{\xi_s}(d\omega)ds \\
   \le &2K \delta_2\sup_{x\in B(\xi_{h-},\delta_1/2)}\Big[-\log \P_{\delta_x}\left(\int_0^{\delta_2}\langle \textbf{1}_{B^c},\omega_{u}\rangle\,du=0\right)\Big]<\infty.
\end{align*}
Combining \eqref{w5.2}, \eqref{w5.4} and \eqref{w5.5}, we get \eqref{step1}.

{\bf Step 2} Next we  prove that given $\xi^h$,
\begin{equation}\label{step2}
\lim_{t\uparrow h}\frac{\langle f_h, X_t^{2,h,\P}\rangle}{\|\Lambda_t^{h}\|}=0,\quad \textbf{P}_\mu^{(h)}\mbox{-a.s.}\end{equation}
Note that given $\xi^h$,
$$
\langle f_h,X^{2,h,\P}_t\rangle:=\int_{0}^t\int_{\mathbb{D}}\langle f_h,\omega_{t-s}\rangle \mathcal{N}^{2,h}(ds,d\omega),
$$
where $\mathcal{N}^{2,h}(ds,d\omega)$ is a Poisson random measure on $[0,h)\times \mathbb{D}$ with intensity measure
$$
 \textbf{1}_{[0,h)}(s) \textbf{1}_{H(\omega)<h-s}\int_0^\infty y n(\xi_s,dy)\P_{y\delta_{\xi_s}}(X\in d\omega)ds.
 $$
Let $I_2$ be the support of the measure $\mathcal{N}^{2,h}$.
Note that $I_2$ is a random countable subset of $[0,h)\times \mathbb{D}$.
Using arguments similar to those leading to \eqref{7.11}, we get that
\begin{align*}
  \langle f_h,X_t^{2,h,\P}\rangle
  \le&\epsilon \langle 1,\Lambda^h_t\rangle+2\|f\|_\infty\int_{0}^{h-\delta_2}\int_{\mathbb{D}}\langle 1,\omega_{t-s}\rangle \mathcal{N}^{2,h}(ds,d\omega)\nonumber\\
   &\quad+2\|f\|_\infty\int_{h-\delta_2}^t\int_{\mathbb{D}}\langle \textbf{1}_{\bar{B}^c},\omega_{t-s}\rangle \mathcal{N}^{2,h}(ds,d\omega)\nonumber\\
   =&\epsilon \langle 1,\Lambda^h_t\rangle+2\|f\|_\infty\sum_{(s,\omega)\in (I_2\cap S_1)}\langle 1,\omega_{t-s}\rangle+2\|f\|_\infty\sum_{(s,\omega)\in (I_2\cap S_2)}\langle \textbf{1}_{\bar{B}^c},\omega_{t-s}\rangle\nonumber\\
   =:&\epsilon \langle 1,\Lambda^h_t\rangle+2\|f\|_\infty J_3(t)+2\|f\|_\infty J_4(t),\label{7.13}
\end{align*}
where $S_1$ and $S_2$ are the set defined in \eqref{e:S_1} and \eqref{S_2}.
It follows that
\begin{eqnarray}\label{w5.2'}
  \frac{|\langle f_h,X^{2,h,\P}_t\rangle|}{\|\Lambda^h_t\|} &\le & \epsilon+2\|f\|_\infty\frac{J_3(t)}{\|\Lambda^h_t\|} +2\|f\|_\infty\frac{J_4(t)}{\|\Lambda^h_t\|}.
\end{eqnarray}
So, to prove \eqref{step2}, we only need to prove that
\begin{equation}\label{w5.8}
  \lim_{t\uparrow h}\frac{J_3(t)}{\|\Lambda^h_t\|} =0, \quad \textbf{P}_\mu^{(h)}\mbox{-a.s.},
\end{equation}
and
\begin{equation}\label{w5.9}
  \lim_{t\uparrow h}\frac{J_4(t)}{\|\Lambda^h_t\|} =0,\quad \textbf{P}_\mu^{(h)}\mbox{-a.s.}
\end{equation}
Note that
\begin{align}\label{w5.7}
  &\int_{S_1}\textbf{1}_{[0,h)}(s) \textbf{1}_{H(\omega)<h-s}\int_0^\infty y n(\xi_s,dy)\P_{y\delta_{\xi_s}}(X\in d\omega)ds\nonumber\\
  \le&\int_{0}^{h-\delta_2}\int_0^\infty y n(\xi_s,dy)\P_{y\delta_{\xi_s}}(X_{\delta_2/2}\neq 0)ds\nonumber\\
 \le &\int_{0}^{h-\delta_2}v(\delta_2/2,\xi_s) \int_0^1 y^2n(\xi_s,dy)ds+\int_{0}^{h-\delta_2}\int_1^\infty y n(\xi_s,dy)\,ds\nonumber\\
  \le &Kh (\|v_{\delta_2/2}\|_\infty+1),
\end{align}
where in the second inequality we used the fact that
$$
\P_{y\delta_{\xi_s}}(X_{\delta_2/2}\neq 0)=1-\P_{y\delta_{\xi_s}}(X_{\delta_2/2}= 0)=1-e^{-yv(\delta_2/2,\xi_s)}\le yv(\delta_2/2,\xi_s).
$$
Thus,
$\mathcal{N}^{2,h}(S_1)<\infty$, a.s.,
which implies that \eqref{w5.8}.

To prove \eqref{w5.9} we only need to show that, given $\xi^h$,
\begin{equation}\label{7.24}
 \int_{S_2}\int_0^\infty y n(\xi_s,dy)\P_{y\delta_{\xi_s}}(X\in d\omega)ds<\infty.
\end{equation}
In fact,
\begin{align*}
  &\int_{S_2}\int_0^\infty y n(\xi_s,dy)\P_{y\delta_{\xi_s}}(X\in d\omega)ds \\
  \le  & \int_{h-\delta_2}^h\int_0^\infty y n(\xi_s,dy)\P_{y\delta_{\xi_s}}\Big(\int_0^{\delta_2}\langle \textbf{1}_{B^c},X_{u}\rangle\,du>0\Big)ds\\
   \le &\int_{h-\delta_2}^h\int_1^\infty y n(\xi_s,dy)ds+\int_{h-\delta_2}^h \Big(-\log \P_{\delta_{\xi_s}}\Big(\int_0^{\delta_2}\langle \textbf{1}_{B^c},X_{u}\rangle\,du=0\Big)\Big)\int_0^1 y^2 n(\xi_s,dy)ds\\
   \le &Kh+Kh\sup_{x\in B(\xi_{h-},\delta_1/2)}\Big[-\log \P_{\delta_{x}}\Big(\int_0^{\delta_2}\langle \textbf{1}_{B^c},X_{u}\rangle\,du=0\Big)\Big]<\infty,
\end{align*}
where in the second inequality, we used the fact that
\begin{align*}
  \P_{y\delta_{\xi_s}}\Big(\int_0^{\delta_2}\langle \textbf{1}_{B^c},X_{u}\rangle\,du>0\Big)=&1-\exp\Big\{y \log\P_{\delta_{\xi_s}}\Big(\int_0^{\delta_2}\langle \textbf{1}_{B^c},\omega_{u}\rangle\,du=0\Big)\Big\}\\
  \le& -y\log\P_{\delta_{\xi_s}}\Big(\int_0^{\delta_2}\langle \textbf{1}_{B^c},X_{u}\rangle\,du=0\Big).
\end{align*}
The proof is now complete.

\hfill$\Box$

\bigskip

\noindent\textbf{Proof of Theorem \ref{mthem}:}\quad
Since $\{X_t,t\ge0\}$ is a Hunt process, $t\to X_t$ is right continuous, which implies that
$$\left\{\lim_{t\uparrow H}\frac{X_t}{\|X_t\|}\quad \mbox{ exists}\right\}=\left\{\lim_{t\in\Q \uparrow H}\frac{X_t}{\|X_t\|} \quad\mbox{ exists}\right\},$$
where $\Q$ is the set of all rational numbers in $[0,\infty)$.
And, note that
$$H=\inf\{t\in \Q: \|X_t\|=0\}.$$
Thus, by Corollary \ref{cor2} and Lemma \ref{the:4.1}, we get that
$$
\P_{\mu}\Big[\lim_{t\in\Q\uparrow H}\frac{X_t}{\|X_t\|} \mbox{ exists}\Big]=\int_0^\infty \textbf{P}_{\mu}^{(h)}\Big[\lim_{t\in\Q\uparrow h}\frac{\Lambda^h_t}{\|\Lambda^h_t\|}\mbox{ exists}\Big]\,F_H(dh)=1.
$$
Let $V:=\lim_{t\uparrow H}\frac{X_t}{\|X_t\|}$. Then, for any $f\in\mathcal{B}_b^+(E)$,
by Lemma \ref{the:4.1},
\begin{align*}
  \P_{\mu}[\exp\{-\langle f,V\rangle\}] = & \P_{\mu}\Big[\lim_{t\in\Q\uparrow H}\exp\Big\{-\frac{\langle f,X_t\rangle}{\|X_t\|}\Big\}\Big]\\
   = &\int_0^\infty\lim_{t\in\Q\uparrow h}\textbf{P}_{\mu}^{(h)}\left[\exp\Big(-\frac{\langle f,\Lambda_t^{h}\rangle}{\|\Lambda_t^{h}\|}\Big)\right]\,F_H(dh)\\
   =&\int_0^\infty\Pi_{\nu}^{h}[\exp(-f(\xi_{h-}))]\,F_H(dh).
\end{align*}
Thus, $V$ is a Dirac measure of the form $V=\delta_Z$ and  the law of $Z$ satisfies \eqref{Z:law}.
The proof is now complete.

\hfill$\Box$

\section{Examples}\label{sec:ex}
In this section, we will list some examples that satisfy
 Assumptions {\bf{(H1)}} and {\bf{(H2)}}.
 The purpose of these examples is to show that Assumptions
 {\bf{(H1)}} and {\bf{(H2)}} are satisfied in a lot of cases.
We will not try to give the most general examples possible.

\begin{exam}\label{exam2}
Suppose that $P_t$ is conservative and preserves $C_b(E)$. Let $\mathcal{A}$ be the infinitesimal generator of $P_t$ in $C_b(E)$ and $\mathcal{D}(\mathcal{A})$ be the domain of $\cA$.
Also assume that
$$
\Psi(x,z)=-\alpha(x)z+b(x)z^2,
$$
where $\sup_{x\in E}\alpha(x)\le 0$ and $\inf_{x\in E} b(x)>0$ and $1/b\in \mathcal{D}(\mathcal{A}).$
Then by Remark \ref{psi-domi}, we know that Auumption {\bf{(H1)}} is satisfied.
One can check that
$$
\Big(\frac{b^{-1}(\xi_t)}{b^{-1}(x)}e^{-\int_0^t (b(\xi_s)\cA(1/b)(\xi_s))\,ds},t\ge0\Big)
$$
is a positive martingale under $\Pi_x$.
Thus we define another probability measure $\Pi_x^{1/b}$ by
$$
\frac{\Pi_x^{1/b}}{\Pi_x}\Big|_{\mathcal{F}_t}=\frac{b^{-1}(\xi_t)}{b^{-1}(x)}e^{-\int_0^t (b(\xi_s)\cA(1/b)(\xi_s))\,ds},\quad  t\ge 0.
$$
Let $\cA^{1/b}$ be the infinitesimal generator of $\xi$ under $\Pi^{1/b}$. If
$-\alpha(x)-b(x)\mathcal{A}(1/b)(x)\in \mathcal{D}(\cA^{1/b})$,  then it
follows from \cite[(3.10) and Lemma 4.9]{DH} that $w(t,x)$ exists and satisfies
$$
w(t,x)\le \frac{1}{\inf_{x\in E}b(x)}e^{ct}\frac{\beta_0^2e^{\beta_0t}}{(e^{\beta_0t}-1)^2},
$$
where $c, \beta_0$ are positive constants. Using this, one can check that
Assumption {\bf{(H2)}} is satisfied.
 This example shows that our result covers  Delmas and H\'{e}nard \cite[Corollary 4.14]{DH}.
\end{exam}

Now we give some examples of superprocesses, with general branching mechanisms,
satisfying Assumptions
 {\bf{(H1)}} and {\bf{(H2)}}.

Recall that the general form of branching mechanism is given by
$$
\Psi(x,z)=-\alpha(x)z+b(x)z^2+\int_0^\infty(e^{-yz}-1+yz)n(x,dy).
$$
By \eqref{n:condition}, there exists $K>0$, such that
$$
|\alpha(x)|+b(x)+\int_0^\infty (y\wedge y^2)n(x,dy)\le K.
$$
Thus we have
\begin{equation}\label{e:7.19}
  |\Psi(x,z)|\le 3K(z+z^2), \qquad x\in\R^d.
\end{equation}

In the next two examples, we  always  assume that $E=\R^d$ and that
$\Psi$ satisfies \eqref{com1}  and the following  condition:
for any $M>0$, there exist $c>0$ and $\gamma_0\in (0, 1]$ such that
\begin{equation}\label{local-lip}
|\Psi(x, z)-\Psi(y, z)|\le c|x-y|^{\gamma_0},\quad x,y\in \R^d, \,z\in[0,M].
\end{equation}
By Remark \ref{psi-domi}, condition \eqref{com1} implies that Assumption {\bf{(H1)}} is satisfied. Therefore, in the following examples, we only need to check that   Assumption {\bf{(H2)}} is satisfied.

\begin{exam}\label{exam3}
Assume that the spatial motion $\xi$ is a diffusion on $\R^d$ satisfying the conditions in Example \ref{exam6}. The branching mechanism $\Psi$ is of the form in \eqref{e:branm} and satisfies \eqref{com1} and \eqref{local-lip}. Then the $(\xi,\Psi)$-superprocess $X$ satisfies Assumptions
{\bf{(H1)}} and {\bf{(H2)}}.
\end{exam}

We now proceed to prove the second assertion of the example above.

\begin{lemma}\label{l:lemm1-d}
For $f\in \mathcal{B}_b(\R^d)$ and $x\in \R^d$, the function $t\to P_tf(x)$ is differentiable on $(0, \infty)$. Furthermore,
 there exists a constant $c$ such that for any $t\in(0,1]$, $x\in\R^d$ and $f\in\mathcal{B}_b(\R^d)$,
\begin{equation}\label{e:7.12-d}
  \left|\frac{\partial}{\partial t}P_tf(x)\right|\le c\|f\|_\infty t^{-1}.
\end{equation}
\end{lemma}

\textbf{Proof:}
For $t\in(n,n+1]$, $P_tf(x)=P_{t-n}(P_nf)(x)$.
Thus, we only need to prove the differentiability for $t\in(0,1]$.
It follows from \cite[IV.(13.1)]{LSU} that
\begin{equation}\label{e:7.1-d}
  \left|\frac{\partial}{\partial t}p(t,x, y)\right|\le c_1t^{-\frac{d}2-1}e^{-\frac{c_2|x-y|^2}{t}}.
\end{equation}
Thus by the dominated convergence theorem we have that for all $t\in (0, 1]$
and $x\in \R^d$,
$$
\frac{\partial}{\partial t}P_tf(x)=\int_{\R^d}\frac{\partial}{\partial t}p(t,x,y)f(y)\,dy,$$
and that  for all $t\in (0, 1]$, $x\in \R^d$ and bounded Borel function $f$ on $\R^d$,
$$
  \left|\frac{\partial}{\partial t}P_tf(x)\right|
   \le  c_3 \|f\|_\infty t^{-1}.
$$
The proof is now complete.
\hfill$\Box$

\begin{lemma}\label{l:lemma6-d}
Assume that $f_s(x)$ is uniformly bounded in $(s,x)\in [0,1]\times \R^d$, that is, there is a constant $L>0$ so that, for all $s\in [0, 1]$ and $x\in\R^d$, $|f_s(x)|\le L$. Then there is a constant $c$ such that for any $t\in(0,1]$ and $x,x'\in\R^d$,
$$
\left|\int_0^t P_{t-s}f_s(x)\,ds-\int_0^t P_{t-s}f_s(x')\,ds\right|\le c L(|x-x'|\wedge 1).
$$
\end{lemma}
\textbf{Proof:}
It follows from \cite[IV.(13.1)]{LSU} that
there exist constants $c_1, c_2>0$ such that for all
$t\in (0, 1]$ and $x, x'\in \R^d$,
\begin{equation}\label{e:gradest}
|\nabla_x p(t, x, y)|\le c_1t^{-\frac{d+1}2}e^{-\frac{c_2|x-y|^2}{t}}.
\end{equation}
Thus
\begin{equation}\label{e:lipschitz-d}
\left|p(t,x, y)-p(t,x', y)\right|\le c_3((t^{-1/2}|x-x'|)\wedge 1)
t^{-d/2}\left(e^{-\frac{c_4|x-y|^2}{t}}+e^{-\frac{c_4|x'-y|^2}{t}}\right)\, .
\end{equation}
Hence for any $t\in(0,1]$ and $x,x'\in\R^d$,
\begin{equation}\label{e:holder2-1-d}
  \left|\int_0^t P_{t-s}f_s(x)\,ds-\int_0^t P_{t-s}f_s(x')\,ds\right|
  \le c_5L\int_0^1s^{-1/2}|x-x'|\,ds\le c_6L|x-x'|.
\end{equation}
\hfill$\Box$

\begin{lemma}\label{l:lemma4-d}
Assume that  $f_s(x)$ satisfies the following conditions:
\begin{enumerate}
  \item [(i)] There is a constant $L$ so that, for all $(s,x)\in [0,1]\times \R^d$, $|f_s(x)|\le L$.
  \item [(ii)] For any $t_0\in [0, 1]$, $\lim_{s\to t_0}\sup_{x\in\R^d}|f_s(x)-f_{t_0}(x)|=0$.
  \item [(iii)]
  There exist constants $s_0\in(0, 1)$, $C>0$ and $\gamma\in(0,1]$
  such that for all $s\in [0, s_0]$ and $x, x'\in \R^d$ with $|x-x'|\le 1$,
    \begin{equation}\label{e:17.6-d}
     |f_s(x)-f_s(x')|\le C|x-x'|^{\gamma}.
    \end{equation}
\end{enumerate}
Then, $t\rightarrow \int_0^t P_{t-s}f_s(x)\,ds$ is differentiable on $(0,s_0)$, and for $t\in [0, s_0)$,
\begin{equation}\label{e:7.10-d}
  \frac{\partial}{\partial t}\int_0^t P_{t-s}f_s(x)\,ds=\int_0^t \frac{\partial}{\partial t}P_{t-s}f_s(x)\,ds+f_t(x).
\end{equation}
\end{lemma}
\textbf{Proof:}
Let $G(t,x):=\int_0^t P_{t-s}f_s(x)\,ds$.
First, we will show that for any $x\in \R^d$,
\begin{equation}\label{e:7.2-d}
  \lim_{t\downarrow0}t^{-1}\int_0^t P_{t-s}f_s(x)\,ds=f_0(x).
\end{equation}
Since $f_0\in C_b(\R^d)$,  we have $\lim_{s\to0}P_sf_0(x)=f_0(x)$, which implies that
$$
\lim_{t\to0}t^{-1}\int_0^t P_{t-s}f_0(x)\,ds=\lim_{t\to0}t^{-1}\int_0^t P_{s}f_0(x)\,ds=f_0(x).
$$
Thus, it suffices to prove that
\begin{equation}\label{e:7.3-d}
  \lim_{t\to0}t^{-1}\int_0^t P_{t-s}(f_s-f_0)(x)\,ds=0.
\end{equation}
Notice that
$$
t^{-1}\int_0^t |P_{t-s}(f_s-f_0)(x)|\,ds\le \sup_{s\le t}\|f_s-f_0\|_\infty\to 0,
$$
as $t\to0$.
Thus, \eqref{e:7.2-d} is valid.

For any $0<t<t+r<s_0$, by the definition of $G(t,x)$,
\begin{eqnarray*}
  &&\frac{1}{r}\left(G(t+r,x)-G(t,x)\right)\\
  &=& \frac{1}{r}\int_0^t \Big(P_{t+r-s}f_s(x)-P_{t-s}f_s(x)\Big)\,ds+\frac{1}{r}\int_t^{t+r} P_{t+r-s}f_s(x)\,ds \\
&=& \int_0^t \frac{P_{t+r-s}f_s(x)-P_{t-s}f_s(x)}{r}\,ds+\frac{1}{r}\int_0^{r} P_{r-s}f_{t+s}(x)\,ds \\
&:=&(I)+(II).
\end{eqnarray*}
By \eqref{e:7.2-d}, we have
\begin{equation}\label{e:7.4-d}
  \lim_{r\downarrow0}(II)=f_t(x).
\end{equation}
Now we deal with part $(I)$.
For $0<t<t+r<s_0$, using \eqref{e:7.1}, we obtain that
\begin{eqnarray}\label{e:8.5-d}
  &&\left|\frac{P_{t+r-s}f_s(x)-P_{t-s}f_s(x)}{r}\right|\nonumber\\
  &=&\left|\int_{\R^d}\frac{p(t+r-s,x,y)-p(t-s,x,y)}{r}(f_s(y)-f_s(x))\,dy\right| \nonumber\\
   &\le& c_3\int_{\R^d}|f_s(y)-f_s(x)|(t-s)^{-\frac{d}2-1}e^{-\frac{c_4|x-y|^2}{t-s}}
   \,dy \nonumber\\
   &\le&c_5\int_{\R^d}|x-y|^{\gamma} (t-s)^{-\frac{d}2-1}e^{-\frac{c_4|x-y|^2}{t-s}}\,dy \nonumber\\
   &\le &c_6(t-s)^{\gamma/2-1}.
\end{eqnarray}
Thus, using the dominated convergence theorem, we get that, for any $0\le t<t+r<s_0$,
\begin{equation}\label{e:7.8-d}
  \lim_{r\downarrow0}(I)=\int_0^t \lim_{r\downarrow0}\frac{P_{t+r-s}f_s(x)-P_{t-s}f_s(x)}{r}\,ds=\int_0^t \frac{\partial}{\partial t}P_{t-s}f_s(x)\,ds.
\end{equation}
Combining \eqref{e:7.4-d} and \eqref{e:7.8-d}, we get that
$$
\lim_{r\downarrow 0}\frac{G(t+r,x)-G(t,x)}{r}=\int_0^t \frac{\partial}{\partial t}P_{t-s}f_s(x)\,ds+f_t(x).
$$
Using similar arguments, we can also  show that
$$
\lim_{r\downarrow 0}\frac{G(t,x)-G(t-r,x)}{r}=\int_0^t \frac{\partial}{\partial t}P_{t-s}f_s(x)\,ds+f_t(x).
$$
Thus, \eqref{e:7.10-d} follows immediately. The proof is now complete.

\hfill $\Box$

Recall that $v(s,\cdot)$ is a bounded function and
$$
v(t+s,x)+\int_0^t P_{t-u}(\Psi_{s+u})(x)\,du=P_t(v_s)(x),
$$
where
\begin{equation}\label{psi-u}
\Psi_u(x)=\Psi(x,v(u,x)).
\end{equation}

\begin{lemma}\label{l:lemma7-d}
For any $s>0$, there is a constant $c(s)$ such that for $t\in[0,1/2)$ and $x,y\in\R^d$,
$$
|v_{t+s}(x)-v_{t+s}(y)|\le c(s)|x-y|.
$$
Moreover, $c(s)$ is decreasing in $s>0$.
\end{lemma}
\textbf{Proof:}
Let $e(s):=\frac{1\wedge s}{2}$. Note that $t+e(s)\in(e(s),1)$.
Thus
$$
v(t+s,x)+\int_0^{t+e(s)} P_{t+e(s)-u}(\Psi(\cdot,v_{s-e(s)+u}(\cdot)))(x)\,du=P_{t+e(s)}(v_{s-e(s)})(x).
$$
It follows from \eqref{e:lipschitz-d} that there exists a constant $c_1$ such that
for all $x,y\in\R^d$,
\begin{align}\label{e:7.14-d}
  &\displaystyle|P_{t+e(s)}(v_{s-e(s)})(x)-P_{t+e(s)}(v_{s-e(s)})(y)|\nonumber\\
  \le &\displaystyle c\|v_{s-e(s)}\|_\infty ((t+e(s))^{-1/2}|x-y|)\wedge 1)\nonumber\\
  \le &\displaystyle c\|v_{s-e(s)}\|_\infty (t+e(s))^{-1/2}|x-y|\nonumber\\
  \le &\displaystyle c\|v_{s/2}\|_\infty (e(s))^{-1/2}|x-y|.
\end{align}
Since $v(s-e(s)+u,x)\le v(s-e(s),x)\le v(s/2,x)$,
we have for $u>0$,
$$
\|\Psi(\cdot,v_{s-e(s)+u}(\cdot))\|_\infty\le 3K(\|v_{s/2}\|_\infty+\|v_{s/2}\|_\infty^2).
$$
Applying Lemma \ref{l:lemma6-d}, we get that there is a constant $c_2>0$ such that for $t\in[0,1/2)$ and $x,y\in\R^d$,
\begin{align}\label{e:7.15-d}
  &\left|\int_0^{t+e(s)} P_{t+e(s)-u}(\Psi(\cdot,v_{s-e(s)+u}(\cdot)))(x)\,du-\int_0^{t+e(s)} P_{t+e(s)-u}(\Psi(\cdot,v_{s-e(s)+u}(\cdot)))(y)\,du\right|\nonumber\\
  \le& c_2 3K(\|v_{s/2}\|_\infty+\|v_{s/2}\|_\infty^2)(|x-y|\wedge1).
\end{align}
The conclusions of the lemma now follow immediately from \eqref{e:7.14-d} and \eqref{e:7.15-d}.
\hfill$\Box$

\begin{lemma}\label{l:lemma3-d}
The function $\Psi_u(x)$ given by \eqref{psi-u} satisfies the following two properties:
\begin{enumerate}
  \item [(1)] For any $u_0>0$,
$$
\lim_{u\to u_0}\sup_{x\in\R^d}|\Psi_u(x)-\Psi_{u_0}(x)|=0;
$$
  \item [(2)] For $t_0\in (0, 1)$,  there exists a constant $c>0$ such that
  for any $|x-x'|\le 1$, $s>t_0$ and $t\in [0,1/2]$,
$$
|\Psi_{s+t}(x)-\Psi_{s+t}(x')|\le c|x-x'|^{\gamma_0}.
$$
\end{enumerate}
\end{lemma}
\textbf{Proof:}
(1)
For $z_1<z_2\in[0,a]$, we can easily check that
\begin{align}\label{e:7.16-d}
  &|\Psi(x,z_1)-\Psi(x,z_2)|\nonumber\\
  \le& |\alpha(x)||z_1-z_2|+b(x)|z_1^2-z_2^2|+\int_0^\infty \left| e^{-yz_1}+yz_1-e^{-yz_2}-yz_2\right|n(x,dy)\nonumber\\
   \le&  K(1+2a) |z_1-z_2|+\int_0^\infty (2\wedge(ya))y|z_1-z_2|n(x,dy)\le K(3+3a)|z_1-z_2|,
\end{align}
where in the second inequality above we use the fact that
$$
|\frac{d}{dx}(e^{-x}+x)|=1-e^{-x}\le 2\wedge x.
$$
Thus, for $|u-u_0|\le u_0/2$, we have that
\begin{equation}\label{e:7.17-d}
  |\Psi_u(x)-\Psi_{u_0}(x)|\le 3K(1+\|v_{u_0/2}\|_\infty)|v_{u}(x)-v_{u_0}(x)|.
\end{equation}
Thus, it suffices to show that $t\mapsto v_t(x)$ is continuous on $(0,\infty)$ uniformly in $x$.

It follows from Lemma \ref{l:lemma7-d} that, for any $t>0$, $x\mapsto v_t(x)$ is  uniformly continuous, thus
$$\lim_{r\downarrow0}\|P_rv_t-v_t\|_{\infty}=0.$$
For $r>0$ and $t>0$, we have that
\begin{eqnarray*}
  |v_t(x)-v_{t+r}(x)| &\le& |P_rv_t(x)-v_t(x)|+|\int_0^r P_{r-u}(\Psi_{t+u})(x)\,du| \\
   &\le&   \|P_rv_t-v_t\|_\infty+3K(\|v_t\|_\infty+\|v_t\|_\infty^2)r\to 0,\quad r\downarrow0,
\end{eqnarray*}
where in the last inequality we used  \eqref{e:7.19} and the fact that $v_{t+u}(x)\le v_t(x)$.

The proof of $\lim_{r\downarrow0}\|v_t-v_{t-r}\|_\infty =0$ is similar and omitted. The proof of  part (1) is now complete.

(2) For any $s>t_0$, and $t\in[0,1/2]$, $v(t+s,x)\le \|v_{t_0}\|_\infty$. By our assumption on $\Psi$,
there exist $c_1>0$ and $\gamma_0\in(0,1]$
such that for $|x-y|\le 1$, $s>t_0$ and $t\in[0,1/2]$,
$$
|\Psi(x,v_{s+t}(x))-\Psi(y,v_{s+t}(x))|\le c_1|x-y|^{\gamma_0}.
$$
By Lemma \ref{l:lemma7-d}, there exists $c_2=c_2(t_0)$ such that for $s>t_0$ and $t\in[0,1/2]$,
$$
|v_{s+t}(x)-v_{s+t}(y)|\le c_2|x-y|.
$$
Thus, for $|x-y|\le 1$, $s>t_0$, and $t\in[0,1/2]$,
\begin{eqnarray}\label{e:7.18-d}
 \displaystyle |\Psi_{s+t}(x)-\Psi_{s+t}(y)|&\le& \displaystyle|\Psi(x,v_{s+t}(x))-\Psi(y,v_{s+t}(x))|+|\Psi(y,v_{s+t}(x))-\Psi(y,v_{s+t}(y))| \nonumber\\
   &\le& \displaystyle  |\Psi(x,v_{s+t}(x))-\Psi(y,v_{s+t}(x))|+3K(1+\|v_{t_0}\|_\infty)|v_{s+t}(x)-v_{s+t}(y)|\nonumber\\
   &\le& \displaystyle c_1|x-y|^{\gamma_0}+3K(1+\|v_{t_0}\|_\infty)c_2|x-y|\nonumber\\
   &\le&\displaystyle c_3 |x-y|^{\gamma_0}.
\end{eqnarray}
The proof of (2) is now complete. \hfill$\Box$

\begin{lemma}\label{l:lemma8-d}
The function $t\to v_t(x)$ is differentiable in $(0, \infty)$, and for any $s>0$ and $t\in[0,1/2)$,
 $w(t+s,x)=-\frac{\partial}{\partial t}v_{t+s}(x)$ satisfies that
\begin{equation}\label{e:7.20-d}
  w(t+s,x)
  =-\frac{\partial}{\partial t}P_t(v_s)(x)+\int_0^t \frac{\partial}{\partial t}P_{t-u}(\Psi_{s+u})(x)\,du+\Psi_{t+s}(x).
\end{equation}
Moreover, $t\to w(t,x)$ is continuous and for any $s_0>0$, $\sup_{s>s_0}\sup_{x\in\R^d}w(s,x)<\infty$.
\end{lemma}
\textbf{Proof:}
For any $t,s>0$,
$$
v(t+s,x)+\int_0^t P_{t-u}(\Psi_{s+u})(x)\,du=P_t(v_s)(x).
$$
Thus, combining  Lemmas \ref{l:lemm1-d}, \ref{l:lemma4-d} and  \ref{l:lemma3-d}, \eqref{e:7.20-d} follows immediately.

For fixed $t\in(0,1/2)$, we  deal with the three parts on right hand side of \eqref{e:7.20-d} separately.

Since $t\to v(t,x)$ is continuous, the function $s\to \Psi_{t+s}(x)=\Psi(x,v(t+s,x))$ is continuous and, by \eqref{e:7.19},
\begin{equation}\label{e:10.3-d}
  \sup_{s>t_0}|\Psi_{t+s}(x)|\le 3K(\|v_{t_0}\|_\infty+\|v_{t_0}\|_\infty^2)<\infty.
\end{equation}
By \eqref{e:7.12},
\begin{equation}\label{e:10.2-d}
  \sup_{s>t_0}|\frac{\partial}{\partial t}P_t(v_s)(x)|\le c_4\|v_{t_0}\|_\infty t^{-1}<\infty.
\end{equation}

By \eqref{e:8.5-d} and Lemma \ref{l:lemma3-d} (2), we get that, for any $s>t_0$,
\begin{equation}\label{e:10.4-d}
 \sup_{s>t_0}\sup_{x\in\R^d}|\int_0^t \frac{\partial}{\partial t}P_{t-u}(\Psi_{s+u})(x)\,du|<\infty.
\end{equation}

 Combining \eqref{e:10.3-d} -\eqref{e:10.4-d}, we get that, for $t_0>0$,
$$
\sup_{s>t_0}\sup_{x\in\R^d}w(t+s,x)<\infty,
$$
which implies that, for any $s_0>0$, $\sup_{s>s_0}\sup_{x\in\R^d}w(s,x)<\infty$.
\hfill $\Box$

\bigskip

Now  we  give an  example of a superprocess with discontinuous spatial motion
and general branching mechanism such that Assumptions
{\bf{(H1)}} and {\bf{(H2)}} are satisfied.

\begin{exam}\label{ex:supersbm}
Suppose that $B=\{B_t\}$ is a Brownian motion in $\R^d$ and $S=\{S_t\}$ is an
independent subordinator with Laplace exponent $\varphi$, that is
$$
\mathbb{E} e^{-\lambda S_t}=e^{-t\varphi(\lambda)}, \qquad t>0,\, \lambda>0.
$$
The process $\xi_t=B_{S_t}$ is called a subordinate Brownian motion in $\R^d$.
Subordinate Brownian motions form a large class of L\'evy processes.
When $S$ is an $(\alpha/2)$-stable subordinator, that is, $\varphi(\lambda)=\lambda^{\alpha/2}$, $\xi$ is a symmetric $\alpha$-stable process in $\R^d$.
Suppose that $\Psi$ is of the form in \eqref{e:branm} satisfying
\eqref{com1} and \eqref{local-lip}.
Suppose further that $\varphi$ satisfies the following conditions:
\begin{enumerate}
\item $\int^1_0\frac{\varphi(r^2)}rdr<\infty.$
\item There exist constants $\delta\in (0, 2]$ and $a_1\in (0, 1)$ such that
\begin{equation}\label{e:lowerscaling}
a_1\lambda^{\delta/2}\varphi(r)\le \varphi(\lambda r), \qquad \lambda\ge 1, r\ge 1.
\end{equation}
\end{enumerate}
then $X$ satisfies Assumptions {\bf{(H1)}} and  {\bf{(H2)}}.
\end{exam}

Now we proceed to prove the second assertion of the example above.
The arguments are similar to that for the second assertion of Example \ref{exam3}.
Without loss of generality, we will assume that $\varphi(1)=1$. First we introduce some notation. Put $\Phi(r)=\varphi(r^2)$ and let $\Phi^{-1}$ be the inverse function of $\Phi$. For $t>0$ and $x\in \R^d$, we define
$$
\rho(t, x):={\Phi\left(\left(\frac{1}{\Phi^{-1}(t^{-1})}+|x|\right)^{-1}\right)}\left(\frac{1}{\Phi^{-1}(t^{-1})}+|x|\right)^{-d}\, .
$$
For $t>0$, $x\in \R^d$ and $\beta, \gamma\in \R$, we define
$$
\rho_{\gamma}^{\beta}(t,x):=\Phi^{-1}(t^{-1})^{-\gamma}(|x|^{\beta}\wedge 1) \rho(t,x)\, ,\quad t>0, x\in \R^d\, .
$$
Let $p(t, x, y)=p(t, x-y)$ be the transition density of $\xi$ and let $\{P_t: t\ge 0\}$
be the transition semigroup of $\xi$. It is well known that $\{P_t: t\ge 0\}$
satisfies the strong Feller property, that is, for any $t>0$,  $P_t$ maps bounded
Borel functions on $\R^d$ to bounded continuous functions on $\R^d$.

Now we list some other properties of the semigroup $\{P_t: t\ge 0\}$ which will
be used later.

\begin{lemma}\label{l:lemm1}
For $f\in \mathcal{B}_b(\R^d)$ and $x\in \R^d$, the function $t\to P_tf(x)$ is differentiable on $(0, \infty)$. Furthermore,
 there exists a constant $c$ such that for any $t\in(0,1]$, $x\in\R^d$ and $f\in\mathcal{B}_b(\R^d)$,
\begin{equation}\label{e:7.12}
  \left|\frac{\partial}{\partial t}P_tf(x)\right|\le c\|f\|_\infty t^{-1}.
\end{equation}
\end{lemma}

\textbf{Proof:}
For $t\in(n,n+1]$, $P_tf(x)=P_{t-n}(P_nf)(x)$.
Thus, we only need to prove the differentiability for $t\in(0,1]$.
It follows from \cite[Lemma 3.1(a) and Theorem 3.4]{KSV16} that
\begin{equation}\label{e:7.1}
  \left|\frac{\partial}{\partial t}p(t,x)\right|\le c_1\rho(t, x).
\end{equation}
By \cite[Lemma 2.6(a)]{KSV16}, we have
\begin{equation}\label{e:lemma2.6.a}
\int_{\R^d}\rho(t, x)dx<c_2t^{-1}, \qquad t\in (0, 1].
\end{equation}
Thus by the dominated convergence theorem we have that for all $t\in (0, 1]$
and $x\in \R^d$,
$$
\frac{\partial}{\partial t}P_tf(x)=\int_{\R^d}\frac{\partial}{\partial t}p(t,x,y)f(y)\,dy,$$
and that  for all $t\in (0, 1]$, $x\in \R^d$ and bounded Borel function $f$ on $\R^d$,
$$
  \left|\frac{\partial}{\partial t}P_tf(x)\right|
   \le  c_3 \|f\|_\infty t^{-1}.
$$
The proof is now complete.
\hfill$\Box$

\begin{lemma}\label{l:lemma6}
Assume that $f_s(x)$ is uniformly bounded in $(s,x)\in [0,1]\times \R^d$, that is, there is a constant $L>0$ so that, for all $s\in [0, 1]$ and $x\in\R^d$, $|f_s(x)|\le L$. Then there is a constant $c$ such that for any $t\in(0,1]$ and $x,x'\in\R^d$,
$$
\left|\int_0^t P_{t-s}f_s(x)\,ds-\int_0^t P_{t-s}f_s(x')\,ds\right|\le c L(|x-x'|^{\delta/2}\wedge 1).
$$
\end{lemma}
\textbf{Proof:}
It follows from \cite[Proposition 3.3]{KSV16} that there exists a constant $c_1>0$ such that for all $t\in (0, 1]$ and $x, x'\in \R^d$,
\begin{equation}\label{e:lipschitz}
\left|p(t,x)-p(t,x')\right|
\le c_1\big((\Phi^{-1}(t^{-1})|x-x'|)\wedge 1\big)t\left(\rho(t,x)+\rho(t,x')\right)\, .
\end{equation}
Thus using \eqref{e:lemma2.6.a} we get that for any $t\in(0,1]$ and $x,x'\in\R^d$,
\begin{equation}\label{e:holder2-1}
  \left|\int_0^t P_{t-s}f_s(x)\,ds-\int_0^t P_{t-s}f_s(x')\,ds\right|
  \le c_2L\int_0^t\big((\Phi^{-1}(s^{-1})|x-x'|)\wedge 1\big)\,ds.
\end{equation}
When $|x-x'|<1$,  $\Phi(|x-x'|^{-1})\ge \Phi(1)=1$.
Thus,
$$
\int_0^t\big((\Phi^{-1}(s^{-1})|x-x'|)\wedge 1\big)\,ds\le |x-x'|
\int^1_{(\Phi(|x-x'|^{-1}))^{-1}}\Phi^{-1}(s^{-1})ds+\big(\Phi(|x-x'|^{-1})\big)^{-1}.
$$
It is well known that $\varphi$,  the Laplace exponent of a subordinator,
satisfies
$$
\varphi(\lambda r)\le \lambda\varphi(r), \qquad \lambda\ge 1, r>0.
$$
Using this, we immediately get that
$$
\Phi^{-1}(\lambda r)\ge \lambda^{1/2}\Phi^{-1}(r), \qquad \lambda\ge 1, r>0.
$$
For $s\in [(\Phi(|x-x'|^{-1})^{-1}, 1]$, by taking $r=s^{-1}$ and $\lambda=s\Phi(|x-x'|^{-1})$ in the display above, we get
$$
\Phi^{-1}(s^{-1})
\le |x-x'|^{-1}s^{-1/2}(\Phi(|x-x'|^{-1}))^{-1/2}.
$$
Therefore
\begin{align*}
 &|x-x'|\int^1_{(\Phi(|x-x'|^{-1}))^{-1}}\Phi^{-1}(s^{-1})ds\\
& \le (\Phi(|x-x'|^{-1}))^{-1/2}\int^1_{(\Phi(|x-x'|^{-1}))^{-1}}s^{-1/2}ds
 \le c_3(\Phi(|x-x'|^{-1}))^{-1/2}.
\end{align*}
Consequently for all $t\in (0, 1]$ and $x, x'\in \R^d$ with $|x-x'|<1$, we have
$$
\int_0^t((\Phi^{-1}(t^{-1})|x-x'|)\wedge 1)\,ds\le c_4(\Phi(|x-x'|^{-1}))^{-1/2}.
$$
By taking $r=1$ and $\lambda=|x-x'|^{-1}$ in \eqref{e:lowerscaling}, we get
$$
a_1|x-x'|^{-\delta}\le \Phi(|x-x'|^{-1}).
$$
Thus for all $t\in (0, 1]$ and $x, x'\in \R^d$ with $|x-x'|<1$, we have
\begin{equation}\label{e:holder2-2}
\int_0^t((\Phi^{-1}(s^{-1})|x-x'|)\wedge 1)\,ds\le c_4a_1^{-1/2}|x-x'|^{\delta/2}.
\end{equation}
Combining \eqref{e:holder2-1} and \eqref{e:holder2-2}, we immediately get the desired conclusion.
\hfill$\Box$

\begin{lemma}\label{l:lemma4}
Assume that  $f_s(x)$ satisfies the following conditions:
\begin{enumerate}
  \item [(i)] There is a constant $L$ so that, for all $(s,x)\in [0,1]\times \R^d$, $|f_s(x)|\le L$.
  \item [(ii)] For any $t_0\in [0, 1]$, $\lim_{s\to t_0}\sup_{x\in\R^d}|f_s(x)-f_{t_0}(x)|=0$.
  \item [(iii)] There exist constants $s_0\in(0, 1)$, $\gamma\in (0,\delta/2]$ and $C>0$ such that for all $s\in [0, s_0]$ and $x, x'\in \R^d$ with $|x-x'|\le 1$,
    \begin{equation}\label{e:17.6}
     |f_s(x)-f_s(x')|\le C|x-x'|^{\gamma}.
    \end{equation}
\end{enumerate}
Then, $t\rightarrow \int_0^t P_{t-s}f_s(x)\,ds$ is differentiable on $(0,s_0)$, and for $0\le t<s_0$,
\begin{equation}\label{e:7.10}
  \frac{\partial}{\partial t}\int_0^t P_{t-s}f_s(x)\,ds=\int_0^t \frac{\partial}{\partial t}P_{t-s}f_s(x)\,ds+f_t(x).
\end{equation}
\end{lemma}
\textbf{Proof:}
Let $G(t,x):=\int_0^t P_{t-s}f_s(x)\,ds$.
For any $0<t<t+r<s_0$, by the definition of $G(t,x)$,
\begin{eqnarray*}
  &&\frac{1}{r}\left(G(t+r,x)-G(t,x)\right)\\
  &=& \frac{1}{r}\int_0^t \Big(P_{t+r-s}f_s(x)-P_{t-s}f_s(x)\Big)\,ds+\frac{1}{r}\int_t^{t+r} P_{t+r-s}f_s(x)\,ds \\
&=& \int_0^t \frac{P_{t+r-s}f_s(x)-P_{t-s}f_s(x)}{r}\,ds+\frac{1}{r}\int_0^{r} P_{r-s}f_{t+s}(x)\,ds \\
&:=&(I)+(II).
\end{eqnarray*}
Using the same arguments as those leading to \eqref{e:7.2-d}, we get
$$  \lim_{t\downarrow0}t^{-1}\int_0^t P_{t-s}f_s(x)\,ds=f_0(x),$$
which implies that
\begin{equation}\label{e:7.4}
  \lim_{r\downarrow0}(II)=f_t(x).
\end{equation}
Now we deal with part $(I)$.
For $0<t<t+r<s_0$, using \eqref{e:7.1}, we obtain that
\begin{eqnarray}\label{e:8.5}
  &&\left|\frac{P_{t+r-s}f_s(x)-P_{t-s}f_s(x)}{r}\right|\nonumber\\
  &=&\left|\int_{\R^d}\frac{p(t+r-s,x,y)-p(t-s,x,y)}{r}(f_s(y)-f_s(x))\,dy\right| \nonumber\\
   &\le& c_3\int_{\R^d}|f_s(y)-f_s(x)|\rho(t-s, x-y)\,dy \nonumber\\
   &\le&c_4\int_{\R^d}\rho^\gamma_0(t-s, x-y)\,dy \nonumber\\
   &\le &c_5(t-s)^{-1}\Phi^{-1}((t-s)^{-1})^{-\gamma},
\end{eqnarray}
where in the last inequality we used \cite[Lemma 2.6(a)]{KSV16}. It follows from
\cite[Lemma 2.3]{KSV16} that
\begin{equation}\label{e:rsnew}
\int^t_0(t-s)^{-1}\Phi^{-1}((t-s)^{-1})^{-\gamma}ds\le c_6\Phi^{-1}(t^{-1})^{-\gamma}.
\end{equation}
Thus, using the dominated convergence theorem, we get that, for any $0\le t<t+r<s_0$,
\begin{equation}\label{e:7.8}
  \lim_{r\downarrow0}(I)=\int_0^t \lim_{r\downarrow0}\frac{P_{t+r-s}f_s(x)-P_{t-s}f_s(x)}{r}\,ds=\int_0^t \frac{\partial}{\partial t}P_{t-s}f_s(x)\,ds.
\end{equation}
Combining \eqref{e:7.4} and \eqref{e:7.8}, we get that
$$
\lim_{r\downarrow 0}\frac{G(t+r,x)-G(t,x)}{r}=\int_0^t \frac{\partial}{\partial t}P_{t-s}f_s(x)\,ds+f_t(x).
$$
Using similar arguments, we can also  show that
$$
\lim_{r\downarrow 0}\frac{G(t,x)-G(t-r,x)}{r}=\int_0^t \frac{\partial}{\partial t}P_{t-s}f_s(x)\,ds+f_t(x).
$$
Thus, \eqref{e:7.10} follows immediately. The proof is now complete.\hfill $\Box$

\bigskip

\begin{lemma}\label{l:lemma7}
For any $s>0$, there is a constant $c(s)$ such that for $t\in[0,1/2)$ and $x,y\in\R^d$,
$$
|v_{t+s}(x)-v_{t+s}(y)|\le c(s)|x-y|^{\delta/2}.
$$
Moreover, $c(s)$ is decreasing in $s>0$.
\end{lemma}
\textbf{Proof:}
The proof of this lemma is similar as that of Lemma \ref{l:lemma7-d}. We use Lemma \ref{l:lemma6} instead of Lemma \ref{l:lemma6-d}.
Here we omit the details.
\hfill$\Box$

\begin{lemma}\label{l:lemma3}
The function $\Psi_u(x)$ satisfies the following two properties:
\begin{enumerate}
  \item [(1)] For any $u_0>0$,
$$
\lim_{u\to u_0}\sup_{x\in\R^d}|\Psi_u(x)-\Psi_{u_0}(x)|=0;
$$
  \item [(2)] For $t_0\in (0, 1)$,  there exists a constant $c>0$ and $\gamma_1\in(0,\delta/2]$ such that
  for any $|x-x'|\le 1$, $s>t_0$ and $t\in [0,1/2]$,
$$
|\Psi_{s+t}(x)-\Psi_{s+t}(x')|\le c|x-x'|^{\gamma_1}.
$$
\end{enumerate}
\end{lemma}
\textbf{Proof:}
The proof of part (1) is exactly the same as that of part (1)
of Lemma \ref{l:lemma3-d}.

Using arguments similar to that in the proof of  part (2) of Lemma \ref{l:lemma3-d}
and using Lemma \ref{l:lemma7} instead of Lemma \ref{l:lemma7-d}, we can get the result in part (2).
Here we omit the details.
\hfill$\Box$

\begin{lemma}\label{l:lemma8}
The functiom $t\to v_t(x)$ is differentiable in $(0, \infty)$, and for any $s>0$ and $t\in[0,1/2)$,
$w(t+s,x)=-\frac{\partial}{\partial t}v_{t+s}(x)$ satisfies that
\begin{equation}\label{e:7.20}
  w(t+s,x)
  =-\frac{\partial}{\partial t}P_t(v_s)(x)+\int_0^t \frac{\partial}{\partial t}P_{t-u}(\Psi_{s+u})(x)\,du+\Psi_{t+s}(x).
\end{equation}
Moreover, $t\to w(t,x)$ is continuous and for any $s_0>0$, $\sup_{s>s_0}\sup_{x\in\R^d}w(s,x)<\infty$.
\end{lemma}
\textbf{Proof:}
Combining  Lemmas \ref{l:lemm1}, \ref{l:lemma4} and  \ref{l:lemma3},
and using arguments similar to that in the proof of
Lemma \ref{l:lemma8-d}, Lemma \ref{l:lemma8} follows immediately.
\hfill $\Box$

\begin{remark}
Actually, by the same arguments and the results from \cite{KSV16}, one check
that in the example above, we could have replaced the subordinate Brownian motion
by the non-symmetric jump process considered there, which contains the non-symmetric stable-like process discussed in \cite{ChZh}.
\end{remark}

Let $L$ be as in Example \ref{exam6}.
Let $E$ be a bounded smooth domain in $\R^d$ and let $p(t, x, y)$ be the
Dirichlet heat kernel of $L$ in $E$. It follows from \cite[Theorem 2.1, p. 247]{GM} that there exist $c_i>0, i=1, 2, 3, 4$, such that for all $t\in (0, 1]$,
\begin{align*}
& \left|\frac{\partial}{\partial t}p(t,x, y)\right|\le c_1t^{-\frac{d}2-1}e^{-\frac{c_2|x-y|^2}{t}},\quad\mbox{ and }\\
&|\nabla_x p(t, x, y)|\le c_3t^{-\frac{d+1}2}e^{-\frac{c_4|x-y|^2}{t}}.
\end{align*}
Using these instead of \eqref{e:7.1-d} and \eqref{e:gradest}, and repeating the arguments for Example \ref{exam3}, we can get the following example.

\begin{exam}
Assume that $E$ be is bounded smooth domain in $\R^d$
and that the spatial motion is $\xi^E$,
 which is the diffusion $\xi$ of Example \ref{exam6} killed upon exiting $E$.
The branching mechanism $\Psi$ is of the form in \eqref{e:branm} and satisfies \eqref{com1} and \eqref{local-lip} on $E$. Then the $(\xi^E,\Psi)$-superprocess $X$ satisfies Assumptions
{\bf{(H1)}} and {\bf{(H2)}}.
\end{exam}

\begin{singlespace}

\end{singlespace}

\vskip 0.2truein
\vskip 0.2truein

\noindent{\bf Yan-Xia Ren:} LMAM School of Mathematical Sciences \& Center for
Statistical Science, Peking
University,  Beijing, 100871, P.R. China. Email: {\texttt
yxren@math.pku.edu.cn}

\smallskip
\noindent {\bf Renming Song:} Department of Mathematics,
University of Illinois,
Urbana, IL 61801, U.S.A.
Email: {\texttt rsong@math.uiuc.edu}

\smallskip

\noindent{\bf Rui Zhang:}  School of Mathematical Sciences, Capital Normal
University,  Beijing, 100048, P.R. China. Email: {\texttt
zhangrui27@cnu.edu.cn}

\end{document}